\documentclass{article}

\usepackage[
  left=1.25in,
  right=1.25in,
  top=1.2in,
  bottom=1.2in
]{geometry}

\usepackage{tikz}

\usepackage{amsmath, amssymb, amsthm}
\usepackage{mathtools}
\usepackage{xcolor}
\usepackage[normalem]{ulem}
\usepackage{hyperref}
\usepackage{enumitem}
\usepackage{booktabs}
\usepackage{array}
\usepackage[style=numeric-comp,hyperref=true,doi=false,url=false,
            isbn=false,giveninits=true,sorting=none,
            block=none,backend=bibtex,maxnames=99]{biblatex}
\newcommand{\N}{\mathbb{N}}
\newcommand{\R}{\mathbb{R}}
\newcommand{\C}{\mathbb{C}}
\newcommand{\Z}{\mathbb{Z}}
\newcommand{\T}{\mathbb{T}}
\newcommand{\frakg}{\mathfrak{g}}
\newcommand{\frakk}{\mathfrak{k}}
\newcommand{\fraks}{\mathfrak{s}}
\newcommand{\frakz}{\mathfrak{z}}
\newcommand{\frakt}{\mathfrak{t}}
\newcommand{\frakh}{\mathfrak{h}}
\newcommand{\frakn}{\mathfrak{n}}
\newcommand{\frakw}{\mathfrak{w}}
\newcommand{\frakc}{\mathfrak{c}}
\newcommand{\xalg}{\mathfrak X_{\mathrm{alg}}}
\newcommand{\xalgz}{\mathfrak X^{0}_{\mathrm{alg}}}
\newcommand{\xtrig}{\mathfrak X_{\mathrm{trig}}}
\newcommand{\sX}{\mathsf{X}}
\newcommand{\Ocal}{\mathcal{O}}
\newcommand{\cL}{\mathcal{L}}
\newcommand{\Sym}{\operatorname{Sym}}
\newcommand{\xc}[1]{\vspace{.1cm}
\noindent {\em #1} }

\newtheorem{theorem}{Theorem}
\newtheorem{proposition}[theorem]{Proposition}
\newtheorem{corollary}[theorem]{Corollary}
\newtheorem{lemma}{Lemma}
\theoremstyle{definition}
\newtheorem{remark}{Remark}
\newtheorem{example}{Example}
\newtheorem{definition}{Definition}

\theoremstyle{plain}

\newtheorem{mainthm}{Theorem}

\newtheorem{maincor}{Corollary}

\newtheorem{maincorC}{Corollary}

\title{Algebraic Vector Fields on Compact Homogeneous Spaces are $2$-generated}

\author{Yiyang Jiang$^*$ \quad and \quad Xudong Chen\footnote{Y. Jiang and X. Chen are with the Department of Electrical \& Systems Engineering, Washington University, St.~Louis, MO 63130, USA. Emails: {\tt\small \char123 j.yiyang, cxudong\char125 @wustl.edu}.
}}

\date{}

\begin{document}

\maketitle

{\renewcommand{\thefootnote}{}}

\begin{abstract}
Let $M$ be a compact, connected manifold on which a compact Lie group $K$ acts smoothly and transitively.
We show that the Lie algebra of algebraic vector fields on $M$, i.e., of smooth vector fields whose translates under $K$ span finite-dimensional spaces, is generated by two vector fields.
Consequently, the subgroup generated by the flows of these two vector fields is dense in the identity component of the diffeomorphism group of~$M$.
For spheres, explicit generators are obtained from earlier results on spherical vector fields.
\end{abstract}

\section{Introduction and main results}\label{sec:introduction}

Let $M$ be a compact, connected manifold of positive dimension, and $K$ be a compact Lie group, with Lie algebra~$\frakk$, that acts smoothly and transitively on~$M$.
We refer to $M$ as a {\it compact homogeneous space}.
Examples include spheres, real, complex, and quaternionic projective spaces, Grassmannians, flag manifolds, Stiefel manifolds, compact Lie groups, tori, and products such as $S^1\times S^2$.
 For $k\in K$ and a smooth vector field $V$ on $M$, let $k\cdot V$ be the push-forward of $V$ by the diffeomorphism $x\mapsto k\cdot x$; we call the vector fields $k\cdot V$, for $k\in K$, the {\it translates} of $V$ under $K$.
The vector field $V$ is said to be {\it $K$-finite} if the span of its translates under $K$ is finite-dimensional.
We call the $K$-finite vector fields {\it algebraic}, since, for every equivariant embedding of $M$ into a finite-dimensional real $K$-module $\mathbb V$, they are the restrictions to $M$ of the polynomial vector fields on $\mathbb V$ that are tangent to $M$ (see Subsection~\ref{ssec:reductions}).
Let $\xalg(M)$ be the space of algebraic vector fields on $M$, which is a Lie algebra over~$\R$.
Given a set $\mathcal S \subseteq \xalg(M)$, let $\operatorname{Lie}(\mathcal S)$ be the Lie algebra generated by~$\mathcal S$.

The main result of the paper is the following:

\begin{mainthm}\label{thm:main}
For any compact homogeneous space $M$, there exist two algebraic vector fields $W_1$ and $W_2$ such that $\xalg(M)=\operatorname{Lie}\{W_1,W_2\}$.
\end{mainthm}

The two vector fields in Theorem~\ref{thm:main} are constructed as follows.
For $X\in\frakk$, let $R_X$ be the vector field on $M$ given by $R_X(x):=\frac{d}{dt} \big|_{t=0}e^{-tX}\cdot x$, and let $\sX_0:=\{R_X\mid X\in\frakk\}$.
We fix a suitable equivariant embedding of $M$ into a finite-dimensional real $K$-module $\mathbb V$ (see Sections~\ref{sec:degree-growth} and~\ref{sec:two-generators}), and let $\sX_1$ be the span of the vector fields $\phi R_X$, where $\phi$ is the restriction to $M$ of a linear function on~$\mathbb V$ and $X\in\frakk$.
If the semi-simple part $[\frakk,\frakk]$ of $\frakk$ acts nontrivially on $M$, i.e., if $R_X\neq0$ for some $X\in[\frakk,\frakk]$, then we take $W_1=R_T$ and $W_2=R_E+w$.
Here, $T$ belongs to a maximal abelian subalgebra of~$\frakk$, $E$ is the sum of nonzero simple root vectors of $[\frakk,\frakk]$ and their complex conjugates, and $w\in\sX_1$ (see Theorem~\ref{thm:two-generators}); the values of finitely many weights at $T$ are required to be pairwise distinct,
  and $w$ is in general position, i.e., it belongs to a nonempty Zariski open subset of a subspace of~$\sX_1$.
Otherwise, $M$ is a torus on which the identity component of $K$ acts by translations, and explicit generators are given in Proposition~\ref{prop:torus-two}.
In both cases, the components of $W_2$ in the eigenspaces of $\operatorname{ad}_{W_1}$ are obtained by polynomial interpolation, as in the proof of Kuranishi's theorem that every real or complex semi-simple Lie algebra is generated by two elements~\cite{Kuranishi49}, \cite[Theorem~6]{Kuranishi51}.

We now describe a consequence of Theorem~\ref{thm:main} for diffeomorphism groups.
Let $\mathfrak X(M)$ be the space of smooth vector fields on $M$, $\operatorname{Diff}(M)$ be the group of smooth diffeomorphisms of $M$, and $\operatorname{Diff}_0(M)$ be the subgroup of $\operatorname{Diff}(M)$ comprising elements that are smoothly isotopic to the identity map $\operatorname{id}_M$.
We equip $\mathfrak X(M)$ and $\operatorname{Diff}(M)$ with the $\mathrm{C}^\infty$-topology, i.e., the topology of uniform convergence of derivatives of all orders, where vector fields and diffeomorphisms of $M$ are regarded as maps from $M$ to a Euclidean space in which $M$ is embedded, and derivatives are taken in local charts of~$M$; this topology does not depend on the choices made.
The Lie algebra $\xalg(M)$ is dense in $\mathfrak X(M)$; this is a form of the Peter--Weyl theorem (see~\cite[p.~1205]{Luk75Izv}), and we include the argument in the proof of Corollary~\ref{cor:B} below.
 Since every vector field on $M$ is complete, it follows from Theorem~\ref{thm:main} that every compact homogeneous space has the strong density property in the sense of~\cite{CofD3gen} with two vector fields, i.e., there are two complete vector fields on $M$ that generate a dense Lie subalgebra of $\mathfrak X(M)$.
 Let $W_1$ and $W_2$ be algebraic vector fields on $M$ with $\xalg(M)=\operatorname{Lie}\{W_1,W_2\}$, which exist by Theorem~\ref{thm:main}, and let $\mathcal W$ be the subgroup of $\operatorname{Diff}_0(M)$ generated by $e^{tW_i}$, for $t\in\R$ and $i=1,2$.
It is a classical fact that the flows of finitely many vector fields that generate a dense Lie subalgebra of $\mathfrak X(M)$ generate a dense subgroup of $\operatorname{Diff}_0(M)$: 
the vector fields whose flows lie in a closed subgroup of $\operatorname{Diff}(M)$ form a closed Lie subalgebra of $\mathfrak X(M)$~\cite[Theorem~1]{Luk74}, \cite[Corollary~1]{Luk81}, and $\operatorname{Diff}_0(M)$ is generated by flows~\cite{Thurston74}, \cite[Corollary~2]{Luk81}; see also~\cite[Theorem~2.5]{Varolin01}.
We use this fact in the form given in~\cite[Theorem~1 and Remark~1]{AK26}, which builds upon~\cite{AS20},
    and obtain the following corollary of Theorem~\ref{thm:main}:

\begin{maincor}\label{cor:B}
The subgroup $\mathcal W$ is dense in $\operatorname{Diff}_0(M)$ with respect to the $\mathrm{C}^\infty$-topology.
\end{maincor}

\begin{proof}
We first show that $\xalg(M)$ is dense in $\mathfrak X(M)$.
Let $V\in\mathfrak X(M)$, and $dg$ be the normalized Haar measure on $K$.
Since $(g,x)\mapsto(g\cdot V)(x)$ is smooth on $K\times M$, the map $g\mapsto g\cdot V$ from $K$ to $\mathfrak X(M)$ is continuous; hence, $V$ is approximated by the vector fields $\int_K\psi(g)\,g\cdot V\,dg$, for continuous functions $\psi\geq0$ on $K$ with $\int_K\psi\,dg=1$ and sufficiently small support around the identity.
By differentiation under the integral sign, $V_\varphi:=\int_K\varphi(g)\,g\cdot V\,dg$ is a smooth vector field on $M$ for every continuous function $\varphi$ on $K$, and $\varphi\mapsto V_\varphi$ is continuous with respect to uniform convergence on $K$.
The matrix coefficients $g\mapsto\langle\pi(g)a,b\rangle$ of the orthogonal representations $\pi$ of $K$ on Euclidean spaces form a subalgebra of $\mathrm{C}(K)$ that contains the constants and separates points (see, e.g.,~\cite[Chapter~III]{BtD}), so every continuous function on $K$ is a uniform limit of such matrix coefficients by the Stone--Weierstrass theorem.
For a matrix coefficient $\varphi$, the vector field $V_\varphi$ is algebraic: for $h\in K$, we have that $h\cdot V_\varphi=\int_K\varphi(h^{-1}g)\,g\cdot V\,dg$, and the functions $g\mapsto\varphi(h^{-1}g)$, for $h\in K$, span a finite-dimensional space.
Thus, $\xalg(M)$ is dense in $\mathfrak X(M)$, and so is $\operatorname{Lie}\{W_1,W_2\}=\xalg(M)$.

Now, let $f\in\operatorname{Diff}_0(M)$, and $(P_t)_{t\in[0,1]}$ be a smooth isotopy from $\operatorname{id}_M$ to $f$; then, $f$ is the time-one map of the smooth time-dependent vector field $V_t:=\dot P_t\circ P_t^{-1}$.
By~\cite[Theorem~1]{AK26}, the flow of every smooth time-dependent vector field on $[0,1]$ with values in the closure of $\operatorname{Lie}\{W_1,W_2\}$, such as $V_t$, can be approximated, uniformly on $[0,1]$ with respect to the $\mathrm{C}^\infty$-topology on $\operatorname{Diff}(M)$, by the flows of time-dependent vector fields $u_1(t)W_1+u_2(t)W_2$ with measurable, bounded controls $u_1$ and $u_2$.
By~\cite[Remark~1]{AK26}, these controls can be chosen piecewise constant with at most one nonzero component at each time, so that the time-one maps of the corresponding flows are finite compositions of fixed-time maps of the flows of $W_1$ and $W_2$, i.e., elements of $\mathcal W$.
Thus, $f$ belongs to the closure of $\mathcal W$, and since $\mathcal W\subseteq\operatorname{Diff}_0(M)$, the subgroup $\mathcal W$ is dense in $\operatorname{Diff}_0(M)$.
\end{proof}

Corresponding results for semi-simple matrix groups and for $S^2$ have been obtained in~\cite[Section~4.2]{CofD3gen} and~\cite[Section~4]{CofDS2} by the same argument, with~\cite[Theorem~2.5]{Varolin01} and with direct estimates for flows, respectively, in place of~\cite[Theorem~1]{AK26}.

 If $W_1=R_X$ for some $X\in\frakk$, as for the generators constructed in this paper, then the flow of $W_1$ is given by the action of a one-parameter subgroup of the compact Lie group~$K$, and Corollary~\ref{cor:B} can be sharpened as follows; here, a subset of a topological group, or of a topological Lie algebra, is said to {\it topologically generate} it if the subgroup, or the Lie subalgebra, generated by the subset is dense (see~\cite[pp.~1203--1204]{Luk75Izv}).

\begin{maincorC}\label{cor:three}
Suppose that $W_1=R_X$ for some $X\in\frakk$.
Then, there is $a>0$ such that $\operatorname{Diff}_0(M)$ is topologically generated by the three diffeomorphisms $e^{aW_1}$, $e^{W_2}$, and $e^{\sqrt2W_2}$.
\end{maincorC}

\begin{proof}
Since $\Z+\sqrt2\,\Z$ is dense in $\R$ and $t\mapsto e^{tW_2}$ is continuous, the maps $e^{W_2}$ and $e^{\sqrt2 W_2}$ topologically generate the closure of the flow of $W_2$.
Since $W_1=R_X$, its flow is the image of the one-parameter subgroup $\{\exp(tX) \mid t\in\R\}$ of $K$ under the continuous homomorphism from $K$ to $\operatorname{Diff}(M)$ given by the action; hence, its closure $\mathcal H$ is the image of the closure of this one-parameter subgroup, which is a torus, and $\mathcal H$ is a torus as well.
For every nontrivial character $\chi$ of $\mathcal H$, the homomorphism $t\mapsto\chi(e^{tW_1})$ is nontrivial, so $\chi(e^{aW_1})=1$ only for $a$ in a countable set; for $a>0$ outside the union of these countably many sets, no nontrivial character of $\mathcal H$ is trivial at $e^{aW_1}$, and hence, the closed subgroup generated by $e^{aW_1}$ is $\mathcal H$ (see, e.g.,~\cite[Chapter~III]{BtD}).
Consequently, the closure of the subgroup generated by $e^{aW_1}$, $e^{W_2}$, and $e^{\sqrt2 W_2}$ contains the flows of $W_1$ and $W_2$, and the assertion follows from Corollary~\ref{cor:B}, since these three diffeomorphisms belong to $\operatorname{Diff}_0(M)$.
\end{proof}

In the language of geometric control theory~\cite{AS22,AK26}, Corollary~\ref{cor:B} states that the control-linear system $\dot x(t)=u_1(t)W_1(x(t))+u_2(t)W_2(x(t))$ on $M$, with two scalar control inputs $u_1$ and $u_2$, is approximately controllable on $\operatorname{Diff}_0(M)$ (cf.~\cite[Definition~1]{AK26}): by the proof of Corollary~\ref{cor:B}, every element of $\operatorname{Diff}_0(M)$ is a limit, in the $\mathrm{C}^\infty$-topology, of time-one maps of the flows of this system for piecewise constant controls that have at most one nonzero component at each time.
Approximate controllability on $\operatorname{Diff}_0(M)$ in the $\mathrm{C}^0$-topology, for any finite family of vector fields whose Lie algebra approximates every vector field in a suitable sense, has been shown in~\cite[Theorem~5.1]{AS22},
    and its $\mathrm{C}^\infty$ version is~\cite[Proposition~1 and Remark~1]{AK26}.

Finite systems of generators for Lie algebras of vector fields and for diffeomorphism groups of compact manifolds have been studied by Lukatskii~\cite{Luk74,Luk75Dokl,Luk75Izv,Luk77,Luk78,Luk79,Luk81}, mostly in the topological sense.
For a compact homogeneous space $M$ with $K$ connected, it has been shown in~\cite[Theorem~2]{Luk75Izv} that $\mathfrak X(M)$ is topologically generated by the vector fields $R_X$, for $X$ in a set of generators of $\frakk$, together with any $a_M$ analytic vector fields from an everywhere dense subset of $\mathfrak X(M)^{a_M}$, and that $\operatorname{Diff}_0(M)$ is topologically generated by $2(a_M+1)$ diffeomorphisms.
 Here, $a_M$ is the smallest integer such that $c_Q \leq  a_M\dim Q$ for every irreducible complex $K$-module $Q$, where $c_Q$ is the multiplicity of $Q$ in the complexification of $\xalg(M)$; since $\xalg(M)\neq0$ and $c_Q \leq \dim Q\cdot\dim M$ for every $Q$, we have that $1\leq a_M\leq\dim M$~\cite[p.~1204]{Luk75Izv}.
For example, $a_M=1$ if $M$ is a compact, irreducible symmetric space and $K$ is the identity component of its isometry group~\cite[pp.~1206--1207]{Luk75Izv}, and $a_M=\dim M$ if $M$ is a torus acting on itself by translations~\cite[p.~1212]{Luk75Izv}.
Since $a_M\geq1$, the systems of generators in~\cite[Theorem~2]{Luk75Izv} consist of at least four diffeomorphisms and, if $\dim M\geq2$, of at least three vector fields, as $\frakk$ is then not generated by a single element.
Moreover, the $a_M$ analytic vector fields are not all algebraic: they are obtained from~\cite[Theorem~1]{Luk75Izv}, so that the span of their translates under $K$ is dense in $\mathfrak X(M)$, whereas the translates of finitely many algebraic vector fields span a finite-dimensional space.

 In the works cited above, generators of $\xalg(M)$ as a Lie algebra, i.e., without closure, are given only for spheres, real projective spaces, and tori~\cite{Luk74,Luk77,Luk79}.
For the sphere $S^{n-1}$ with $n\geq3$, it follows from~\cite[Theorem~3]{Luk74} and~\cite[Theorem~1 and Lemma~2]{Luk77} that the divergence-free part of $\xalg(S^{n-1})$ and $\xalg(S^{n-1})$ itself are generated by $n$ and by $n+1$ explicit vector fields, respectively; these results are recalled in Section~\ref{sec:spheres}.
Analogous results hold for real projective spaces~\cite[Theorem~1 and Lemma~3]{Luk77}, and for the torus $\T^m$ with $m\geq2$, the Lie algebra $\xalg(\T^m)$ is generated by $m+2$ explicit vector fields~\cite[Theorem~1]{Luk79}.
In these three cases, $\operatorname{Diff}_0(M)$ is topologically generated by two diffeomorphisms~\cite[Theorems~3 and~5]{Luk77}, \cite[Theorem~2]{Luk79}.
For an arbitrary compact manifold $M$, it has been shown in~\cite[Theorem~5]{Luk81}, and announced in~\cite[Theorem~4]{Luk78}, that $\mathfrak X(M)$ and $\operatorname{Diff}_0(M)$ are topologically finitely generated, by a construction that uses an atlas of $M$; see also~\cite[Theorem~1.4]{Akhmedov15} for the groups of $\mathrm{C}^k$-diffeomorphisms, $k\geq1$.

 Theorem~\ref{thm:main} concerns the Lie algebra $\xalg(M)$ itself: for every compact homogeneous space $M$, the two generators $W_1$ and $W_2$ are algebraic, and every algebraic vector field on $M$ is a linear combination of iterated Lie brackets of $W_1$ and $W_2$, so that no closure is involved.
Since $\xalg(M)$ is dense in $\mathfrak X(M)$, it follows that $\mathfrak X(M)$ is topologically generated by two vector fields, and by Corollaries~\ref{cor:B} and~\ref{cor:three}, $\operatorname{Diff}_0(M)$ is topologically generated by the flows of two vector fields, as well as by three diffeomorphisms; the numbers of generators for spheres are listed in Table~\ref{tab:counts}.
The proof of Theorem~\ref{thm:main} is based on the degree-growth theorem of Section~\ref{sec:degree-growth} and on the equivariant embeddings constructed in Section~\ref{sec:two-generators}.

For the sphere $S^{n-1}\subseteq\R^n$ with $K=\mathrm{SO}_n(\R)$, the Lie algebra $\xalg(S^{n-1})$ comprises the polynomial vector fields that are tangent to $S^{n-1}$ (see Section~\ref{sec:spheres}); for $S^2$,
    two explicit generators of $\xalg(S^2)$ have been given in~\cite{CofDS2}, and another proof of this result is given in Section~\ref{sec:spheres}.
For a compact Lie group $K$ acting on itself by left multiplication, $\xalg(K)$ is the Lie algebra of algebraic vector fields on $K$ regarded as a matrix group (see Remark~\ref{rem:groups}). It has been shown in~\cite{CofD3gen} that, for every connected, semi-simple matrix group over $\R$ or $\C$, the Lie algebra of algebraic vector fields is generated by three complete vector fields; for compact, connected, semi-simple groups, Theorem~\ref{thm:main} shows that two vector fields suffice. 
 Further results on finite systems of generators for Lie algebras of vector fields, on $\C^d$, $\mathrm{SL}_2(\C)$, Danielewski surfaces, and smooth complex affine varieties~\cite{AndristCn,BeldievPogudin,AndristSLtwo,AndristNew,AndristDanielewski,Bouaziz26}, and on $\R^d$, tori, and $S^2$~\cite[Section~6]{AS22}, are discussed in~\cite[Section~1]{CofD3gen}.

 Theorem~\ref{thm:main} provides two generators of $\xalg(M)$ on every compact homogeneous space $M$; they are explicit on tori and on spheres (see Proposition~\ref{prop:torus-two} and Section~\ref{sec:spheres}), whereas in general, the second generator $R_E+w$ involves a vector field $w$ in general position, for which we do not know explicit choices.
In Section~\ref{sec:spheres}, we show that, on $S^{n-1}$ with $n\geq3$, two explicit generators of the divergence-free part of $\xalg(S^{n-1})$ and of $\xalg(S^{n-1})$ itself are obtained by combining two results of~\cite{Luk77} with the arguments in the proof of Kuranishi's theorem; for $S^2$,
  one of these results requires an argument that is not covered by the proof in~\cite{Luk77}, which we supply in Subsection~\ref{ssec:S2}.

The remainder of the paper is organized as follows: We gather key notations used throughout the paper at the end of this section.
In Section~\ref{sec:degree-growth}, we collect basic facts on algebraic functions and vector fields on compact homogeneous spaces,
    and we recall, in Lemma~\ref{lem:interpolation}, the interpolation argument in the proof of Kuranishi's theorem.
In the same section, we extend the degree-growth theorem of~\cite{CofD3gen}, which concerns functions on semi-simple Lie groups, to Lie algebras of the form $\Sym(W)\otimes\frakg$, where $\frakg$ is a complex semi-simple Lie algebra and $W$ is a finite-dimensional $\frakg$-module, and we show that $\xalg(M)=\operatorname{Lie}(\sX_0+\sX_1)$ for every equivariant embedding of $M$, with a single exception on the circle.
In Section~\ref{sec:two-generators}, we prove Theorem~\ref{thm:main}.
 In Section~\ref{sec:spheres}, we consider spheres, obtain explicit generators from the results of~\cite{Luk77}, and prove that the conformal algebra of $S^2$ is a maximal Lie subalgebra of $\xalg(S^2)$.

\xc{Notation:} We use $\N$ to denote the set of positive integers and $\N_0:=\N\cup\{0\}$, and $\mathrm{i}$ to denote the imaginary unit.
For a real vector space $U$, let $U^\C:=U \otimes_\R \C$ be its complexification, and $\sigma$ be the complex conjugation of $U^\C$ with respect to~$U$.   
A complex subspace $U'$ of $U^\C$ with $\sigma(U')=U'$ is the complexification of the real subspace $U'\cap U$ of its $\sigma$-fixed points.
For a space $\mathcal F$ of (real- or complex-valued) functions and a space $\mathcal V$ of vector fields on a manifold, we let $\mathcal F\mathcal V$ be the (real or complex) span of the vector fields $fV$ with $f\in\mathcal F$ and $V\in\mathcal V$.  

Throughout the paper, we let $K$ be a compact Lie group with Lie algebra $\frakk$, and $M$ be a compact, connected manifold of positive dimension with a smooth transitive action $K\times M\to M$, $(k,x)\mapsto k\cdot x$.
For $x\in M$, we let $K_x$ be the stabilizer of $x$ and $\frakk_x$ be its Lie algebra.
For $k\in K$, a smooth (real- or complex-valued) function $f$ on $M$, and a smooth vector field $V$ on $M$, we let
$$
    (k\cdot f)(x):= f(k^{-1}\cdot x) \quad\mbox{and}\quad
    (k\cdot V)(k\cdot x):= d_x k\,(V(x)), \quad\mbox{for } x\in M,
$$
where $d_x k$ is the differential at $x$ of the diffeomorphism $y\mapsto k\cdot y$.
Then, $k\cdot(fV)=(k\cdot f)(k\cdot V)$, $(k\cdot V)(k\cdot f)=k\cdot(Vf)$, and $k\cdot[U,V]=[k\cdot U,k\cdot V]$ for smooth vector fields $U$ and $V$ on $M$.
Recall that $R_X(x)=\frac{d}{dt} \big|_{t=0} e^{-tX}\cdot x$ for $X\in\frakk$ and that $\sX_0=\{R_X \mid X\in\frakk\}$.
Since the flow of $R_X$ is given by $x\mapsto e^{-tX}\cdot x$, we have that
\begin{equation}\label{eq:R-derivatives}
R_X f=\frac{d}{dt} \Big|_{t=0} e^{tX}\cdot f \quad\mbox{and}\quad
[R_X,V]=\frac{d}{dt} \Big|_{t=0} e^{tX}\cdot V.
\end{equation}
In particular, $[R_X,R_Y]=R_{[X,Y]}$ for all $X,Y\in\frakk$, and $k\cdot R_X=R_{\operatorname{Ad}_kX}$ for all $k\in K$.

\xc{Algebraic functions and vector fields:}
 As in the case of vector fields, a smooth function $f$ on $M$ is said to be {\it $K$-finite} if the span of its translates $k\cdot f$, for $k\in K$, is finite-dimensional.
Let $\Ocal(M)$ be the space of real-valued $K$-finite smooth functions on $M$, which we call {\it algebraic} functions.
In~\cite[p.~1204]{Luk75Izv}, the $K$-finite smooth sections of homogeneous vector bundles over $M$, such as functions and vector fields, are called spherical.
By the identities above, $\Ocal(M)$ is a unital subalgebra of $\mathrm{C}^\infty(M)$, $\xalg(M)$ is a Lie subalgebra of $\mathfrak X(M)$ and an $\Ocal(M)$-module that contains $\sX_0$, and $Vf\in\Ocal(M)$ for all $V\in\xalg(M)$ and $f\in\Ocal(M)$.
We write $\Ocal^\C(M)$ and $\xalg^\C(M)$ for the complexifications of $\Ocal(M)$ and $\xalg(M)$.
Since $k\cdot \bar f=\overline{k\cdot f}$, the space $\Ocal^\C(M)$ comprises the complex-valued $K$-finite smooth functions on $M$,
on which the elements of $\xalg^\C(M)$ act by derivations.
We let $R_{X+\mathrm{i}Y}:=R_X+\mathrm{i}R_Y\in\xalg^\C(M)$ for $X,Y\in\frakk$, so that $X\mapsto R_X$ is a homomorphism of complex Lie algebras from $\frakk^\C$ to $\xalg^\C(M)$.

\section{Degree growth on compact homogeneous spaces}\label{sec:degree-growth}

In this section, we show that $\xalg(M)=\operatorname{Lie}(\sX_0+\sX_1)$ for every equivariant embedding of~$M$, with a single exception on the circle; see Theorem~\ref{thm:every-embedding} below.

\xc{Equivariant embeddings:}
 Let $\mathbb V$ be a finite-dimensional real $K$-module, i.e., a finite-dimensional real vector space with a continuous linear action of $K$, and let $v\in \mathbb V$ be such that $K_v=K_{x_0}$ for some $x_0\in M$.
Then, the map $k\cdot x_0\mapsto k\cdot v$ from $M$ to $\mathbb V$ is well defined, injective, and $K$-equivariant.
Its differential at $x_0$ is injective, since the differential at the identity of $k\mapsto k\cdot v$ is $X\mapsto X\cdot v$, whose kernel is the Lie algebra of $K_v =K_{x_0}$, i.e., $\frakk_{x_0}$; by equivariance, the map is an immersion, and hence, it is an embedding of the compact manifold $M$ onto the orbit $K\cdot v$.
We call the pair $(\mathbb V,v)$ an {\it equivariant embedding} of~$M$.
 It has been shown in~\cite{Mostow57,Palais57} that every compact manifold with a smooth (not necessarily transitive) action of $K$ admits a $K$-equivariant embedding into a finite-dimensional real $K$-module; the image $v$ of any point $x_0\in M$ under such an embedding of $M$ satisfies $K_v=K_{x_0}$, since the embedding is injective and $K$-equivariant.
Thus, equivariant embeddings of $M$ exist.
For example, $(\R^n,e_n)$ is an equivariant embedding of the sphere $S^{n-1}$ with $K=\mathrm{SO}_n(\R)$.
The Grassmannian of $k$-dimensional subspaces of $\R^n$, $\C^n$, or $\mathbb H^n$ is the orbit of an orthogonal projection of rank $k$ under $\mathrm{SO}_n(\R)$, $\mathrm{SU}_n$, or $\mathrm{Sp}_n$, respectively, acting by conjugation on the real vector space of symmetric, Hermitian, or quaternionic Hermitian matrices.
The flag manifold $K/K_X$, for $X\in\frakk$, is the adjoint orbit $K\cdot X\subseteq\frakk$, so $(\frakk,X)$ is an equivariant embedding of $K/K_X$, where $K$ acts on $\frakk$ by the adjoint action.
The Stiefel manifold of orthonormal $k$-frames in $\R^n$, $k<n$, is the orbit of $(e_1,\ldots,e_k)$ in $(\R^n)^k$ under $\mathrm{SO}_n(\R)$, and the product $S^1\times S^{n-1}$ is the orbit of $(e_2,e_n)$ in $\R^2\oplus\R^n$ under $\mathrm{SO}_2(\R)\times\mathrm{SO}_n(\R)$.
Compact Lie groups and tori are treated in Remark~\ref{rem:groups} and Subsection~\ref{ssec:tori}.

Given an equivariant embedding $(\mathbb V,v)$, we identify $M$ with $K\cdot v$.
Let $\mathcal A_1$ be the space of restrictions to $M$ of linear functions on~$\mathbb V$, and $\mathcal A_0:=\R$ be the space of constant functions.
For each $j\geq2$, let $\mathcal A_j$ be the span of the products $\phi_1\cdots\phi_j$ with $\phi_i\in\mathcal A_1$, and for each $j\geq1$, let $\sX_j:=\mathcal A_j\sX_0$.
Let $\mathbb V^*$ be the dual space of $\mathbb V$, on which $K$ acts by $(k\cdot\ell)(y):=\ell(k^{-1}\cdot y)$.
For $\ell\in \mathbb V^*$ and $X\in\frakk$, we have that $R_X(\ell |_M)=(X\cdot\ell) |_M$, where $(X\cdot\ell)(y):=-\ell(X\cdot y)$.
Thus, $\mathcal A_1$ is invariant under $\sX_0$, and so is each $\mathcal A_j$ by the Leibniz rule.
Also, $\mathcal A_j\subseteq\Ocal(M)$, since the translates of the restriction of a polynomial of degree at most $j$ on $\mathbb V$ are restrictions of polynomials of degree at most $j$.
For smooth (real- or complex-valued) functions $\xi,\eta$ and smooth vector fields $U,U'$ on $M$, we have that
\begin{equation}\label{eq:bracket-degree}
[\xi U,\eta U']=\xi\eta[U,U']+\xi(U\eta)U'-\eta(U'\xi)U.
\end{equation}
In particular, $[\xi R_X,\eta R_Y]\in\sX_{a+b}$ for $\xi\in\mathcal A_a$, $\eta\in\mathcal A_b$, and $X,Y\in\frakk$, so $\sum_{j\geq0}\sX_j$ is a Lie subalgebra of $\xalg(M)$.

For the remainder of the section, we use the above notation for a given equivariant embedding $(\mathbb V,v)$ of~$M$.
The main result of the section is the following:

\begin{theorem}\label{thm:every-embedding}
For every equivariant embedding $(\mathbb V,v)$ of $M$, it holds that $\xalg(M)=\operatorname{Lie}(\sX_0+\sX_1)$, unless $M$ is a circle and $\mathcal A_1\subseteq\operatorname{span}\{1,\cos\theta,\sin\theta\}$ for an angle coordinate $\theta$ on $M$ in which the identity component of $K$ acts by rotations.
In the exceptional case, $\operatorname{Lie}(\sX_0+\sX_1)=\operatorname{span}\{\partial_\theta,\cos\theta\,\partial_\theta,\sin\theta\,\partial_\theta\}$, which is isomorphic to $\mathfrak{sl}_2(\R)$.
\end{theorem}

The exceptional case occurs, for example, for the circle $S^1\subseteq\R^2$ with $K=\mathrm{SO}_2(\R)$ and $(\mathbb V,v)=(\R^2,e_2)$ (see Section~\ref{sec:spheres}); it does not occur if the semi-simple part $[\frakk,\frakk]$ of $\frakk$ acts nontrivially on $M$ (see Theorem~\ref{thm:reductive} and the proof of Theorem~\ref{thm:every-embedding} below).

The proof of Theorem~\ref{thm:every-embedding} occupies the rest of the section.
 In Subsection~\ref{ssec:reductions}, we proceed in three steps.
First, we show that $\xalg(M)=\sum_{j\geq0}\sX_j$, so that the first assertion of Theorem~\ref{thm:every-embedding} amounts to the statement that $\sX_j\subseteq\operatorname{Lie}(\sX_0+\sX_1)$ for every $j\geq2$, apart from the exceptional case.
Second, we show that $K$ can be assumed to be connected.
Third, we distinguish two cases, according to whether the semi-simple part $[\frakk,\frakk]$ of $\frakk$ acts nontrivially on $M$ or not.
The two cases are treated in Subsections~\ref{ssec:reductive} and~\ref{ssec:tori}, respectively.

\subsection{Reductions}\label{ssec:reductions}

\xc{The identity $\xalg(M)= \sum_{j\geq0} \sX_j$:}
We have the following lemma, which identifies $\Ocal(M)$ with the algebra generated by the restrictions of linear functions:

\begin{lemma}\label{lem:kfinite}
For every equivariant embedding $(\mathbb V,v)$ of $M$, it holds that $\Ocal(M)=\sum_{j\geq0} \mathcal A_j$.
\end{lemma}

\begin{proof}
We have shown that $\sum_{j\geq0}\mathcal A_j\subseteq\Ocal(M)$.
For the converse, let $\mathcal P$ be the complex span of $\bigcup_{j\geq0}\mathcal A_j$,
which is a unital subalgebra of the algebra $\mathrm{C}(M,\C)$ of continuous complex-valued functions on $M$.
It is closed under complex conjugation and separates the points of $M\subseteq\mathbb V$, so it is dense in $\mathrm{C}(M,\C)$ with respect to the supremum norm by the Stone--Weierstrass theorem.
Let $dg$ be the normalized Haar measure on $K$.
For an irreducible complex $K$-module $Q$ with character $\chi_Q$, let
$$
P_Qf:=\dim Q \int_K \overline{\chi_Q(g)}\,(g\cdot f)\,dg, \quad\mbox{for } f\in\mathrm{C}(M,\C).
$$
Then, $P_Q$ is a continuous linear operator on $\mathrm{C}(M,\C)$, and on every finite-dimensional $K$-invariant subspace, it is the projection onto the isotypic component of type $Q$ (see, e.g.,~\cite[Chapter~II]{BtD}).
The image $\mathrm{C}(M,\C)_Q$ of $P_Q$ is finite-dimensional.
 To wit, fix a basis of $Q$, and let $q_{ij}(g)$, for $1\leq i,j\leq\dim Q$, be the entries of the matrix of $g\in K$ acting on $Q$, so that $\chi_Q(kg)= \sum_{i,j}q_{ij}(k)\,q_{ji}(g)$ for all $k,g\in K$.
For $f\in\mathrm{C}(M,\C)$ and $k\in K$, we have that $(P_Qf)(k\cdot v)=\dim Q\int_K\overline{\chi_Q(g)}\,f(g^{-1}k\cdot v)\,dg$.
Replacing $g$ by $kg$, which leaves the Haar measure invariant, we obtain that
$$
(P_Qf)(k\cdot v)
=\dim Q\int_K \overline{\chi_Q(kg)}\,f(g^{-1}\cdot v)\,dg
=\dim Q\sum_{i,j} \overline{q_{ij}(k)} \int_K\overline{q_{ji}(g)}\,f(g^{-1}\cdot v)\,dg.
$$
Thus, for every $F\in\mathrm{C}(M,\C)_Q$, the function $k\mapsto F(k\cdot v)$ on $K$ belongs to the span of the $(\dim Q)^2$ functions $\overline{q_{ij}}$; since $M=K\cdot v$, the linear map $F\mapsto(k\mapsto F(k\cdot v))$ is injective, so $\dim\mathrm{C}(M,\C)_Q\leq(\dim Q)^2$.
Since $\mathcal P$ is the union of the finite-dimensional $K$-invariant subspaces $\sum_{j \leq d}\mathcal A_j\otimes\C$,
for $d\geq0$, we have that $P_Q( \mathcal P)\subseteq\mathcal P$.
Since $\mathcal P$ is dense and $P_Q$ is continuous, $P_Q(\mathcal P)$ is dense in $\mathrm{C}(M,\C)_Q$, and hence, $P_Q(\mathcal P)=\mathrm{C}(M,\C)_Q$.
Now, let $f\in\Ocal(M)$.
Then, $f$ is the sum of its projections $P_Qf$ onto the isotypic components of the finite-dimensional $K$-module spanned by the translates of $f$, and each $P_Qf$ belongs to $\mathrm{C}(M,\C)_Q=P_Q(\mathcal P)\subseteq\mathcal P$.
Thus, $f\in\mathcal P$, and since $f$ is real-valued, $f\in\sum_{j\geq0}\mathcal A_j$.
\end{proof}

We next show that the algebraic vector fields are the combinations of the vector fields $R_X$ with coefficients in $\Ocal(M)$:

\begin{lemma}\label{lem:tangent}
It holds that $\xalg(M)=\Ocal(M)\sX_0$.
\end{lemma}

\begin{proof}
Since $k\cdot(fR_X)=(k\cdot f)R_{\operatorname{Ad}_kX}$, we have that $\Ocal(M)\sX_0\subseteq\xalg(M)$.
For the converse, fix an $\operatorname{Ad}(K)$-invariant inner product $\langle\cdot,\cdot\rangle$ on $\frakk$,
and for $x\in M$, let $\frakk_x^\perp$ be the orthogonal complement of $\frakk_x$.
Since the action is transitive, the linear map $\frakk_x^\perp\to T_xM$, $X\mapsto R_X(x)$, is an isomorphism; let $s_x\colon T_xM\to\frakk_x^\perp$ be its inverse.
Then, $s_x$ depends smoothly on $x$, since $x\mapsto\frakk_x^\perp$ is a smooth subbundle of the trivial bundle $M\times\frakk$, namely, the orthogonal complement of the kernel of the surjective smooth bundle map $(x,X)\mapsto R_X(x)$ from $M\times\frakk$ to $TM$.
Since $R_{\operatorname{Ad}_kX}(k\cdot x)=d_xk\,(R_X(x))$ and $\operatorname{Ad}_k\frakk_x^\perp=\frakk_{k\cdot x}^\perp$, we have that $s_{k\cdot x}(d_xk\,(u))=\operatorname{Ad}_k\,s_x(u)$ for all $k\in K$ and $u\in T_xM$.
Let $V\in\xalg(M)$, let $X_1,\ldots,X_d$ be an orthonormal basis of $\frakk$, and let $f_j(x):=\langle s_x(V(x)),X_j\rangle$.
Then, $f_j$ is smooth and $V=\sum_{j=1}^df_jR_{X_j}$, since $V(x)=R_{s_x(V(x))}(x)$.
For $k\in K$, we have that $V(k^{-1}\cdot x)=d_xk^{-1}\,((k \cdot V)(x))$, and hence,
$$
(k\cdot f_j)(x)=\langle s_{k^{-1} \cdot x}(V(k^{-1}\cdot x)),X_j\rangle
=\langle\operatorname{Ad}_{k^{-1}} s_x((k \cdot V)(x)),X_j\rangle
=\langle s_x((k\cdot V)(x)),\operatorname{Ad}_k X_j\rangle.
$$
The right-hand side belongs to the span of the functions $x\mapsto\langle s_x(U(x)),Y\rangle$, for $U$ in the span of $\{k\cdot V\mid k\in K\}$ and $Y\in\frakk$, which is finite-dimensional. Thus, $f_j\in\Ocal(M)$.
\end{proof}

 By Lemmas~\ref{lem:kfinite} and~\ref{lem:tangent}, $\xalg(M)=\sum_{j\geq0}\sX_j$ for every equivariant embedding of $M$.
Thus, for every equivariant embedding $(\mathbb V,v)$ of $M$, the algebraic functions on $M$ are the restrictions to $M$ of the polynomial functions on $\mathbb V$, and the algebraic vector fields on $M$ are the restrictions to $M$ of the polynomial vector fields on $\mathbb V$, i.e., of the maps from $\mathbb V$ to $\mathbb V$ with polynomial components, that are tangent to~$M$.
To wit, for $\phi\in\mathcal A_j$ and $X\in\frakk$, let $p$ be a polynomial on $\mathbb V$ with $p|_M=\phi$; then, $\phi R_X$ is the restriction to $M$ of the polynomial vector field $y\mapsto-p(y)\,X\cdot y$, since $R_X(y)=-X\cdot y$ for $y\in M$, and this polynomial vector field is tangent to $M$.
Conversely, let $P$ be a polynomial vector field on $\mathbb V$ that is tangent to $M$; then, $k\cdot(P|_M)=(k\circ P\circ k^{-1})|_M$ for $k\in K$, and the polynomial vector fields $k\circ P\circ k^{-1}$ have the same degree as $P$, so the translates of $P|_M$ span a finite-dimensional space, i.e., $P|_M\in\xalg(M)$.

\begin{remark}\label{rem:groups}
Let $K$ act on itself by left multiplication, and fix a faithful orthogonal representation of $K$ on a Euclidean space $U$,
    so that $K$ is a matrix group in $\operatorname{End}(U)$ and $R_X(g)=-Xg$ for $X\in\frakk$ and $g\in K$.
Then, $\xalg(K)$ is the Lie algebra of algebraic vector fields on the matrix group $K$ in the sense of~\cite{CofD3gen}, i.e., of vector fields whose matrix entries are polynomials in the matrix entries.
To wit, the pair $(\operatorname{End}(U),\operatorname{id}_U)$, with $K$ acting by left multiplication, is an equivariant embedding of $K$.
By Lemmas~\ref{lem:kfinite} and~\ref{lem:tangent}, $\xalg(K)$ comprises the vector fields $g\mapsto\sum_jf_j(g)X_jg$, where $X_j\in\frakk$ and the functions $f_j$ are polynomials in the matrix entries.
Since $g^{-1}=g^{\mathsf T}$, we have that $gX=(gXg^{\mathsf T})g$ and $Xg=g(g^{\mathsf T}Xg)$ for $X\in\frakk$.
Thus, $\xalg(K)$ is also spanned by the left-invariant vector fields $g\mapsto gX$ multiplied by polynomials in the matrix entries.
Conversely, if a vector field $V$ on $K$ has matrix entries that are polynomials in the matrix entries,
  then so does $g\mapsto g^{\mathsf T}V(g)\in\frakk$, so $V$ is of the above form.
\end{remark}

\xc{Reduction to connected $K$:}
The proof of Theorem~\ref{thm:every-embedding} is based on the fact, recalled below, that if $K$ is connected, then $K$ is the product of a torus contained in its center and its connected subgroup with Lie algebra $[\frakk,\frakk]$.
Since a product of two connected subgroups is connected, such a decomposition requires $K$ to be connected.
We therefore reduce to the case where $K$ is connected.
The Lie algebra $\xalg(M)$ depends on the action of $K$ on $M$; however, it does not change if $K$ is replaced by a closed subgroup that still acts transitively:

\begin{lemma}\label{lem:subgroup}
Let $K'$ be a closed subgroup of $K$ whose action on $M$ is transitive. Then, a smooth function or vector field on $M$ is $K'$-finite if and only if it is $K$-finite.
In particular, $\Ocal(M)$ and $\xalg(M)$ do not change if $K$ is replaced by its identity component, which acts transitively on $M$.
\end{lemma}

\begin{proof}
A $K$-finite function or vector field is clearly $K'$-finite.
Conversely, let $(\mathbb V,v)$ be an equivariant embedding of $M$ for the action of $K$.  Since $K'$ acts transitively,
    $K'\cdot v=K\cdot v$, so $(\mathbb V,v)$, regarded as a $K'$-module, is an equivariant embedding of $M$ for the action of $K'$.
By Lemma~\ref{lem:kfinite}, applied to $K$ and to $K'$, both the $K$-finite and the $K'$-finite smooth functions on $M$ form the algebra $\sum_{j\geq0}\mathcal A_j$.
By Lemma~\ref{lem:tangent}, applied to $K'$, a $K'$-finite vector field is a combination of vector fields $R_X$, with $X$ in the Lie algebra of $K'$, with $K'$-finite, and hence, $K$-finite, coefficients; it is thus $K$-finite.
Finally, let $K_0$ be the identity component of $K$. The orbits of $K_0$ in $M$ are compact, and $M$ is the union of the finitely many orbits $(kK_0)\cdot x$, for $kK_0\in K/K_0$ and a fixed $x\in M$.
Distinct orbits are disjoint and closed, so each orbit is also open, and since $M$ is connected, $K_0$ acts transitively.
\end{proof}

By Lemma~\ref{lem:subgroup}, we can assume that $K$ is connected; this changes neither $\frakk$ nor the spaces $\mathcal A_j$ and $\sX_j$.

\xc{The semi-simple part:}
We now recall the decomposition of $K$ on which the proof is based and fix notation.
For connected $K$, let $\fraks:=[\frakk,\frakk]$ and $\frakz$ be the center of $\frakk$.
Then, $\frakk=\frakz\oplus\fraks$, the Lie algebra $\fraks$ is semi-simple (or zero), and $K=ZS$, where $S$ is the connected subgroup of $K$ with Lie algebra $\fraks$, which is a closed normal subgroup,
and $Z$ is the identity component of the center of $K$, which is a torus (see, e.g.,~\cite[Chapter~V]{BtD}).
We call $S$ the {\it semi-simple part} of $K$, and we let $\sX_0^{\fraks}:=\{R_X\mid X\in\fraks\}$ and $\sX_0^{\fraks,\C}:=\{R_X\mid X\in\fraks^\C\}$.

 The proof of Theorem~\ref{thm:every-embedding} distinguishes two cases, according to whether $S$ acts nontrivially or trivially on $M$.
The following remark, which is not used in the proofs below but explains the need for Subsection~\ref{ssec:algebroid}, describes when $S$ acts transitively on $M$; in this case, $S$ acts nontrivially on $M$, since $M$ has positive dimension.

\begin{remark}\label{rem:transitive}
 Suppose that $K$ is connected.
Then, the group $S$ acts transitively on $M$ if and only if the fundamental group of $M$ is finite.
For simply connected $M$, the transitivity of $S$ is a theorem of Montgomery~\cite{Montgomery50}; we include a proof of the equivalence for completeness.
To wit, fix $x\in M$, and let $L:=SK_x$, which is a closed subgroup of $K$, since $S$ is normal and compact.
The map $kK_x\mapsto kL$ is a fibration $M\to K/L$ with fiber $L/K_x\approx S/(S\cap K_x)$, which is connected.
Since $K=ZS$, the quotient $K/S$ is a torus; since $L\supseteq S$, the subgroup $L$ is normal in $K$, and $K/L$ is a quotient of the torus $K/S$. Thus, $K/L$ is a torus of some dimension $m$.
The exact homotopy sequence of the fibration gives an exact sequence $\pi_1(M)\to\pi_1(K/L)\to\pi_0(L/K_x)=0$, so $\pi_1(M)$ surjects onto $\pi_1(K/L)\approx\Z^m$.
If $\pi_1(M)$ is finite, then $m=0$, so $K=L=SK_x$, and hence, $K\cdot x=S\cdot x$, i.e., $S$ acts transitively.
Conversely, if $S$ acts transitively, then $M\approx S/S_x$, and the exact sequence $\pi_1(S)\to\pi_1(M)\to\pi_0(S_x)$ shows that $\pi_1(M)$ is finite,
    since the fundamental group of the compact semi-simple Lie group $S$ is finite and $S_x$ has finitely many components (see, e.g.,~\cite[Chapter~V]{BtD}).
\end{remark}

For example, the spheres $S^{n-1}$ and the real projective spaces $\R P^{n-1}$ with $n\geq3$, the complex and quaternionic projective spaces, the Grassmannians other than $\R P^1$, the flag manifolds, the Stiefel manifolds of orthonormal $k$-frames, $k<n$, in $\R^n$ with $n\geq3$ and in $\C^n$, and the compact semi-simple Lie groups have finite fundamental groups, so the semi-simple part of $K$ acts transitively on them.
On the other hand, the semi-simple part of $K$ acts nontrivially but not transitively on $S^1\times S^{n-1}$, $n\geq3$, with $K=\mathrm{SO}_2(\R)\times\mathrm{SO}_n(\R)$, and on the unitary group $\mathrm{U}_m$, $m\geq2$, acting on itself by left multiplication; in both cases, the fundamental group is infinite.

The remainder of the section is organized as follows:
In Subsection~\ref{ssec:algebroid}, we extend the degree-growth theorem of~\cite{CofD3gen} to an algebraic setting that does not refer to any manifold; see Theorem~\ref{thm:udg} and Corollary~\ref{cor:udg}.
In Subsection~\ref{ssec:reductive}, we treat the case where $S$ acts nontrivially on $M$ and prove the stronger statement $\xalg(M)=\operatorname{Lie}(\sX_0^{\fraks}+\sX_1)$ of Theorem~\ref{thm:reductive}, which is the form used in Section~\ref{sec:two-generators}.
In Subsection~\ref{ssec:tori}, we treat the case where $S$ acts trivially on $M$, in which $M$ is a torus on which $K$ acts by translations, in Theorem~\ref{thm:torus}; the exceptional case of Theorem~\ref{thm:every-embedding} arises there.
The proof of Theorem~\ref{thm:every-embedding} is given in Subsection~\ref{ssec:proof-every-embedding}.

\subsection{Degree growth in \texorpdfstring{$\Sym(W)\otimes\frakg$}{Sym(W)\textotimes g}}\label{ssec:algebroid}

If $S$ acts transitively on $M$, then~\cite[Theorem~1]{CofD3gen} can be applied on $S$ to the functions $\phi\circ q$, for $\phi\in\mathcal A_1$, where $q\colon S\to M$ is the surjective map $q(g):=g^{-1}\cdot v$; to wit,
    $L_X(\psi\circ q)=(R_X\psi)\circ q$ for every smooth function $\psi$ on $M$ and every $X\in\fraks$, where $L_X$ is the left-invariant vector field on $S$ with $L_X(g)=\frac{d}{dt}\big|_{t=0}ge^{tX}$.
In general, however, the orbits of $S$ in $M$ may be proper submanifolds; by Remark~\ref{rem:transitive},
this happens exactly when the fundamental group of $M$ is infinite.
We therefore use the following algebraic form of~\cite[Theorem~1]{CofD3gen}, which does not refer to any manifold.

Let $\frakg$ be a complex semi-simple Lie algebra, and $W$ be a finite-dimensional $\frakg$-module.
The action of $\frakg$ on $W$ extends uniquely to an action on the symmetric algebra $\Sym(W)$ by derivations, which preserves each homogeneous component $\Sym^n(W)$; we write $X\cdot f$ for the action of $X\in\frakg$ on $f\in\Sym(W)$.
On the vector space $\Sym(W)\otimes\frakg$, we define a bilinear bracket by
\begin{equation}\label{eq:algebroid-bracket}
[f\otimes X,g\otimes Y]:= fg\otimes[X,Y]+f(X\cdot g)\otimes Y-g(Y\cdot f)\otimes X,
\quad\mbox{for } f,g\in\Sym(W)\mbox{ and }X,Y\in\frakg,
\end{equation}
which has the same form as the Lie bracket~\eqref{eq:bracket-degree} of the vector fields $fR_X$ and $gR_Y$.
We have the following fact:

\begin{lemma}\label{lem:algebroid-lie}
The bracket~\eqref{eq:algebroid-bracket} turns $\Sym(W)\otimes\frakg$ into a Lie algebra,
    and $[\Sym^a(W)\otimes \frakg,\Sym^b(W)\otimes \frakg]\subseteq \Sym^{a+b}(W)\otimes \frakg$ for all $a,b\in\N_0$.
\end{lemma}

\begin{proof}
The second assertion holds since $\frakg$ preserves each $\Sym^n(W)$.
For the first assertion, let $\mathfrak D$ be the Lie algebra of derivations of $\Sym(W)$, and consider the Lie algebra $\Sym(W)\otimes\frakg$ with the bracket $[f\otimes X,g\otimes Y]_0:=fg\otimes[X,Y]$.
Every $\delta\in\mathfrak D$ acts on $(\Sym(W)\otimes\frakg,[\cdot,\cdot]_0)$ by the derivation $f\otimes X\mapsto\delta(f)\otimes X$, so the semi-direct product $\mathfrak D\ltimes(\Sym(W)\otimes\frakg)$ is a Lie algebra.
For $f\in\Sym(W)$ and $X\in\frakg$, let $\partial_X\in\mathfrak D$ be the derivation $g\mapsto X\cdot g$ and $f\partial_X\in\mathfrak D$ be the derivation $g\mapsto f(X\cdot g)$,
    and consider the injective linear map $\iota\colon f\otimes X\mapsto(f\partial_X,f\otimes X)$ from $\Sym(W)\otimes\frakg$ to $\mathfrak D\ltimes(\Sym(W)\otimes\frakg)$.
Since $X\mapsto\partial_X$ is a Lie algebra homomorphism, $[f\partial_X,g\partial_Y]=fg\partial_{[X,Y]}+f(X\cdot g)\partial_Y-g(Y\cdot f)\partial_X$; also, the second component of $[\iota(f\otimes X),\iota(g\otimes Y)]$ is $fg\otimes[X,Y]+f(X\cdot g)\otimes Y-g(Y\cdot f)\otimes X$.
Thus, $\iota$ maps the bracket~\eqref{eq:algebroid-bracket} to the Lie bracket of the semi-direct product, which proves the first assertion.
\end{proof}

We identify $\frakg$ with $\Sym^0(W)\otimes\frakg$.
For each $n\in\N$, let $\cL_n$ be the space spanned by the right-normed Lie products $[V_1,[V_2,\cdots,[V_{n-1},V_n]]]$ of depth $n$ with $V_i\in W\otimes\frakg$.
The main result of the subsection is the following:

\begin{theorem}\label{thm:udg}
If $W$ does not contain any nonzero $\frakg$-invariant vector, then for any $n\in\N$,
$$ \cL_n=\Sym^n(W)\otimes\frakg.
$$
In particular, the Lie algebra $\Sym(W)\otimes\frakg$ is generated by $\frakg$ and $W\otimes\frakg$.
\end{theorem}

\begin{proof}

The proof of~\cite[Theorem~1]{CofD3gen} concerns spaces of functions on a semi-simple Lie group, but it uses only the form~\eqref{eq:algebroid-bracket} of the Lie bracket, the action of the Lie algebra on the functions by derivations, and the representation theory of the Lie algebra.
It proceeds by induction on $n$, using that $\cL_n=[W\otimes\frakg,\cL_{n-1}]$, and it shows, for a product $\xi$ of $n$ weight vectors, first that $\xi\otimes H\in\cL_n$ for the elements $H$ of a Cartan subalgebra of $\frakg$ and then that $\xi\otimes X_\alpha\in\cL_n$ for the root vectors $X_\alpha$~\cite[Section~2, Propositions~2--5]{CofD3gen}.
The hypothesis there, that the space of functions does not contain any nonzero constant function, enters only through the fact that the zero weight space is spanned by the vectors $X_\alpha\psi$, for positive roots $\alpha$, root vectors $X_\alpha$, and $\psi$ of weight $-\alpha$.
This fact holds for $W$: by complete reducibility and the hypothesis on $W$, every irreducible summand of $W$ has a nonzero lowest weight, and it is spanned by the iterated images of a lowest weight vector under the root vectors $X_\alpha$ for positive roots $\alpha$, so its zero weight vectors are sums of vectors $X_\alpha\psi$ with $\psi$ of weight $-\alpha$.
Thus, the proof applies with only notational changes, with the space of functions, its homogeneous components,
    and the bracket of vector fields there replaced by $W$, $\Sym^n(W)$, and~\eqref{eq:algebroid-bracket}, respectively.
The last assertion of the theorem follows from the first, since $\Sym(W)\otimes\frakg=\bigoplus_{n\geq0}\Sym^n(W)\otimes\frakg$.
\end{proof}

We now apply the theorem to Lie algebras of vector fields.
Let $N$ be a smooth manifold, $\frakg_0$ be a real semi-simple Lie algebra, and $X\mapsto R_X$ be a Lie algebra homomorphism from $\frakg_0$ into the Lie algebra of smooth vector fields on $N$.
For a complex space $\Psi$ of smooth complex-valued functions on $N$ and $n\in\N$, let $\mathcal A_n(\Psi)$ be the complex span of the products $\phi_1\cdots\phi_n$ with $\phi_i\in\Psi$, and $\mathcal A(\Psi):=\C+\sum_{n\in\N}\mathcal A_n(\Psi)$ be the unital algebra generated by $\Psi$.

\begin{corollary}\label{cor:udg}
Let $\Phi$ be a finite-dimensional complex space of smooth complex-valued functions on $N$ such that $R_X\Phi\subseteq\Phi$ for all $X\in\frakg_0$ and no nonzero $\phi\in\Phi$ satisfies $R_X\phi=0$ for all $X\in\frakg_0$.
Then, for every $n\in\N$, the complex span of the right-normed Lie products of depth $n$ of the vector fields $\phi R_X$, with $\phi\in\Phi$ and $X\in\frakg_0$, is $\mathcal A_n(\Phi)\{R_X\mid X\in\frakg_0\}$.
\end{corollary}

\begin{proof}
Let $\frakg:=\frakg_0^\C$, and extend $X\mapsto R_X$ complex-linearly to a homomorphism into the Lie algebra of complex vector fields, i.e.,
    of complex derivations of $\mathrm{C}^\infty(N,\C)$.
Let $W:=\Phi$, regarded as a $\frakg$-module via $X\cdot\phi:=R_X\phi$; by hypothesis, $W$ does not contain any nonzero $\frakg$-invariant vector.
Let $\varrho\colon\Sym(W)\to\mathrm{C}^\infty(N,\C)$ be the algebra homomorphism that extends the inclusion of $W$.
Then, $\varrho(X\cdot f)=R_X\varrho(f)$ for $f\in\Sym(W)$, since both sides are derivations in $f$ that agree on $W$.
Since~\eqref{eq:bracket-degree} holds on $N$ as well, the linear map $f\otimes X\mapsto\varrho(f)R_X$ is a Lie algebra homomorphism from $\Sym(W)\otimes\frakg$ to the Lie algebra of complex vector fields on $N$.
It maps $\cL_n$ onto the span of the right-normed Lie products of depth $n$ of the vector fields $\phi R_X$, with $\phi\in\Phi$ and $X\in\frakg$,
and it maps $\Sym^n(W)\otimes\frakg$ onto $\mathcal A_n(\Phi)\{R_X\mid X\in\frakg\}$.
Since these complex spans do not change if $X$ ranges over $\frakg_0$ instead of $\frakg$, the corollary follows from Theorem~\ref{thm:udg}.
\end{proof}

For a connected, semi-simple Lie group $N$ with the left-invariant vector fields in place of $R_X$, the hypothesis of Corollary~\ref{cor:udg} on $\Phi$ means that $\Phi$ does not contain any nonzero constant function, and the corollary is the complexified form of~\cite[Theorem~1]{CofD3gen}.
In the corollary, however, the vector fields $R_X$ need not span the tangent spaces of $N$, and the map $\varrho$ in the proof need not be injective; this is the form in which the corollary is applied in the next subsection, where the orbits of $S$ in $M$ may be proper submanifolds.

\subsection{The case where the semi-simple part acts nontrivially}\label{ssec:reductive}

Throughout this subsection, we assume that $K$ is connected and that $S$ acts nontrivially on $M$, and we fix an equivariant embedding $(\mathbb V,v)$ of $M$.
With the notation introduced before Corollary~\ref{cor:udg}, applied to $N=M$, we have that $\mathcal A_n(\mathcal A_1^\C)=\mathcal A_n^\C$ and $\mathcal A(\Psi)\subseteq\Ocal^\C(M)$ for every subspace $\Psi\subseteq\Ocal^\C(M)$.
The main result of the subsection is the following:

\begin{theorem}\label{thm:reductive}
Suppose that $K$ is connected and that $S$ acts nontrivially on $M$.  Then, for every equivariant embedding $(\mathbb V,v)$ of $M$,
$$
  \xalg(M)=\operatorname{Lie}(\sX_0^{\fraks}+\sX_1).
$$
In particular, $\xalg(M)=\operatorname{Lie}(\sX_0+\sX_1)$.
\end{theorem}

The proof occupies the rest of the subsection.
Let
$$
    \cL:=\operatorname{Lie}(\sX_0^{\fraks}+\sX_1) \quad\mbox{and}\quad
    \Phi:=\mathcal A_1^\C,
$$
so that $\cL^\C$ is the complex Lie algebra generated by $\sX_0^{\fraks,\C}+\sX_1^\C$, and $\mathcal A(\Phi)=\Ocal^\C(M)$ by Lemma~\ref{lem:kfinite}.
By Lemma~\ref{lem:tangent}, $\xalg^\C(M)=\Ocal^\C(M)\sX_0^\C$ is spanned by the vector fields $fR_X$, with $f\in\Ocal^\C(M)$ and $X\in\fraks^\C$, together with the vector fields $fR_H$, with $f\in\Ocal^\C(M)$ and $H\in\frakz$.
We show in Lemmas~\ref{lem:semisimple-part} and~\ref{lem:central-part} below that both families of vector fields belong to $\cL^\C$.
We first decompose $\Phi$ according to the action of $S$.

\xc{Decomposition of $\Phi$:}
 For a real $K$-module $U$, the $S$-invariant vectors in $U$ form a $K$-submodule of $U$, since $S$ is normal in $K$.
Thus, if $U$ is irreducible, then either $S$ acts trivially on $U$, or $U$ does not contain any nonzero $S$-invariant vector.
Let $\mathbb V_S$ be the sum of the irreducible $K$-submodules of $\mathbb V$ on which $S$ acts nontrivially, and $\mathbb V_0$ be the space of $S$-invariant vectors in $\mathbb V$.
Both are $K$-submodules of $\mathbb V$, and $\mathbb V=\mathbb V_S\oplus\mathbb V_0$.
To wit, by complete reducibility, $\mathbb V$ is the sum of its irreducible $K$-submodules, each of which is contained in $\mathbb V_S$ or in $\mathbb V_0$; thus, $\mathbb V=\mathbb V_S+\mathbb V_0$.
Moreover, every irreducible $K$-submodule of $\mathbb V_S$ is isomorphic to one of the summands in the definition of $\mathbb V_S$, so $S$ acts nontrivially on it, whereas $S$ acts trivially on every submodule of $\mathbb V_0$; thus, $\mathbb V_S\cap\mathbb V_0$ contains no irreducible $K$-submodule, and hence, it is zero.
Let $\Phi_S$ (resp., $\Phi_0$) be the space of restrictions to $M$ of the complex linear functions on $\mathbb V$ that vanish on $\mathbb V_0$ (resp., on $\mathbb V_S$).
Then, $\Phi=\Phi_S+\Phi_0$, and both spaces are invariant under $\sX_0$.
Since $X\cdot y\in \mathbb V_S$ for $X\in\fraks$ and $y\in \mathbb V$, we have that $R_X\Phi_0=0$ for all $X\in\fraks$.
Also, $\Phi_S$ does not contain any nonzero $\fraks$-invariant function. To wit, $\Phi_S$ is a quotient of the $\fraks$-module of complex linear functions on $\mathbb V_S$, which does not contain any nonzero $\fraks$-invariant element, and by complete reducibility, every quotient of this module is isomorphic to one of its submodules.
Since every finite-dimensional $\fraks$-module $U$ satisfies $U=\fraks\cdot U\oplus U^{\fraks}$, where $U^{\fraks}$ is the space of $\fraks$-invariant elements, it follows that
\begin{equation}\label{eq:s-Phi}
\fraks\cdot\Phi=\fraks\cdot\Phi_S=\Phi_S.
\end{equation}

Since $Z$ is a torus, the operators $R_H$, for $H\in\frakz$, are simultaneously diagonalizable on $\Phi_0$ with purely imaginary eigenvalues.
Let $\eta_1,\ldots,\eta_p$ be a basis of $\Phi_0$ consisting of common eigenvectors of these operators, and let $\psi_j\in(\frakz^\C)^*$ be such that $R_H\eta_j=\psi_j(H)\eta_j$ for all $H\in\frakz$.
For $\kappa=(\kappa_1,\ldots,\kappa_p)\in\N_0^p$, let $\eta^\kappa:=\eta_1^{\kappa_1}\cdots\eta_p^{\kappa_p}$ and $\chi_\kappa:=\sum_{j=1}^p\kappa_j\psi_j$, so that
\begin{equation}\label{eq:eta-weights}
R_X\eta^\kappa=0\quad\mbox{for }X\in\fraks^\C
\quad\mbox{and}\quad
R_H\eta^\kappa=\chi_\kappa(H)\,\eta^\kappa\quad\mbox{for }H\in\frakz.
\end{equation}
If $\chi_\kappa=0$, then $\eta^\kappa$ is annihilated by $\sX_0$, so it is $K$-invariant, and hence, it is a constant function.
Since $\Phi=\Phi_S+\Phi_0$, the algebra $\Ocal^\C(M)=\mathcal A(\Phi)$ is spanned by the products $\eta^\kappa f$ with $\kappa\in\N_0^p$ and $f\in\mathcal A(\Phi_S)$.

The following lemma is where the degree-growth theorem enters:

\begin{lemma}\label{lem:semisimple-part}
It holds that $\eta^\kappa\mathcal A(\Phi_S)\sX_0^{\fraks,\C}\subseteq\cL^\C$ for every $\kappa\in\N_0^p$. Consequently, $\Ocal^\C(M)\sX_0^{\fraks,\C}\subseteq\cL^\C$.
\end{lemma}

\begin{proof}
Let $\cL_S:=\mathcal A(\Phi_S)\sX_0^{\fraks,\C}$, which is a Lie subalgebra of $\xalg^\C(M)$ by~\eqref{eq:bracket-degree}, since $\mathcal A(\Phi_S)$ is invariant under $R_X$ for $X\in\fraks^\C$.
By Corollary~\ref{cor:udg}, applied with $N=M$, $\frakg_0=\fraks$, and $\Phi_S$ in place of $\Phi$, the span of the right-normed Lie products of depth $n$ of the vector fields $\phi R_X$, with $\phi\in\Phi_S$ and $X\in\fraks$, is $\mathcal A_n(\Phi_S)\sX_0^{\fraks,\C}$, for every $n\in\N$.
Since $\cL_S = \sX_0^{\fraks,\C}+\sum_{n\in\N}\mathcal A_n(\Phi_S) \sX_0^{\fraks,\C}$ and the vector fields $\phi R_X$, with $\phi\in\Phi_S$ and $X\in\fraks$, belong to $\sX_1^\C$, it follows that $\cL_S\subseteq\cL^\C$ and that $\cL_S$ is generated as a Lie algebra by $\sX_0^{\fraks,\C}$ and $\Phi_S\sX_0^{\fraks,\C}$.
The second assertion of the lemma follows from the first, since $\Ocal^\C(M)\sX_0^{\fraks,\C}$ is spanned by the vector fields $\eta^\kappa fR_X$ with $\kappa\in\N_0^p$, $f\in\mathcal A(\Phi_S)$, and $X\in\fraks^\C$, which belong to $\eta^\kappa\cL_S$.
We prove the first assertion in two steps.

\xc{Step 1: $\eta^\kappa\sX_0^{\fraks,\C}\subseteq\cL^\C$ for all $\kappa\in\N_0^p$.}
We proceed by induction on $|\kappa|:=\kappa_1+\cdots+\kappa_p$; the case $\kappa=0$ is clear.
If $|\kappa|\geq1$, we write $\eta^\kappa=\eta_j\eta^{\kappa'}$ with $|\kappa'|=|\kappa|-1$.
For $X,Y\in\fraks^\C$, it follows from~\eqref{eq:bracket-degree} and~\eqref{eq:eta-weights} that
$$
[\eta_jR_X,\eta^{\kappa'}R_Y]=\eta^\kappa R_{[X,Y]}.
$$
Since $\eta_jR_X\in\sX_1^\C$ and $\eta^{\kappa'}R_Y\in\cL^\C$ by the induction hypothesis, we have that $\eta^\kappa R_{[X,Y]}\in\cL^\C$; since $[\fraks^\C,\fraks^\C]=\fraks^\C$, the claim follows.

\xc{Step 2: $\eta^\kappa\cL_S\subseteq\cL^\C$ for all $\kappa\in\N_0^p$.}
Fix $\kappa\in\N_0^p$, and let $J_\kappa:=\{A\in\cL_S\mid\eta^\kappa A\in\cL^\C\}$.
For $U\in\cL_S$, we have that $U\eta^\kappa=0$ by~\eqref{eq:eta-weights}, since $U$ is a combination of vector fields $R_X$ with $X\in\fraks^\C$. Thus, for $A\in J_\kappa$ and $U\in\cL_S$,
$$
[\eta^\kappa A,U]=\eta^\kappa[A,U]-(U\eta^\kappa)A=\eta^\kappa[A,U],
$$
which belongs to $\cL^\C$ since $U\in\cL_S\subseteq\cL^\C$. Thus, $J_\kappa$ is an ideal of $\cL_S$.
By Step~1, $J_\kappa$ contains $\sX_0^{\fraks,\C}$, and hence, it contains $[\sX_0^{\fraks,\C},\cL_S]$.
Since $\cL_S$ is generated by $\sX_0^{\fraks,\C}$ and $\Phi_S\sX_0^{\fraks,\C}$, it remains to show that $\Phi_S\sX_0^{\fraks,\C}\subseteq J_\kappa$.

The space $\Phi_S\sX_0^{\fraks,\C}$ is an $\fraks$-module under $X\cdot A:=[R_X,A]$, and it is the direct sum of $\fraks\cdot(\Phi_S\sX_0^{\fraks,\C})=[\sX_0^{\fraks,\C},\Phi_S\sX_0^{\fraks,\C}]\subseteq J_\kappa$ and the space of its elements annihilated by $\fraks$.
Let $\Theta\in\Phi_S\sX_0^{\fraks,\C}$ be annihilated by $\fraks$. It remains to show that $\Theta\in J_\kappa$.
The map $\phi\otimes X\mapsto\phi R_X$ from $\Phi_S\otimes\fraks^\C$ to $\Phi_S\sX_0^{\fraks,\C}$ is a surjective homomorphism of $\fraks$-modules, since $[R_Y,\phi R_X]=(R_Y\phi)R_X+\phi R_{[Y,X]}$.
By complete reducibility, we can write $\Theta=\sum_b\phi_bR_{X_b}$, where $\sum_b\phi_b\otimes X_b\in\Phi_S\otimes\fraks^\C$ is $\fraks$-invariant.
The contraction $\phi\otimes X\mapsto R_X\phi$ from $\Phi_S\otimes\fraks^\C$ to $\Phi_S$ is a homomorphism of $\fraks$-modules, so $\sum_bR_{X_b}\phi_b$ is an $\fraks$-invariant element of $\Phi_S$, which is zero.
Now, let $H\in\frakz$. Since $[R_X,R_H]=0$ for $X\in\fraks$, it follows from~\eqref{eq:bracket-degree} and~\eqref{eq:eta-weights} that
$$
    \sum_b[\eta^\kappa R_{X_b},\phi_bR_H]
    =\eta^\kappa\Bigl(\sum_bR_{X_b}\phi_b\Bigr)R_H-\sum_b\phi_b(R_H\eta^\kappa)R_{X_b}
    =-\chi_\kappa(H)\,\eta^\kappa\Theta,
$$
where $\eta^\kappa R_{X_b}\in\cL^\C$ by Step~1 and $\phi_bR_H\in\sX_1^\C$.
If $\chi_\kappa\neq0$, we choose $H$ with $\chi_\kappa(H)\neq0$ and obtain $\eta^\kappa\Theta\in\cL^\C$.
If $\chi_\kappa=0$, then $\eta^\kappa$ is a constant function, and $\eta^\kappa\Theta\in\cL_S\subseteq\cL^\C$, so $\Theta\in J_\kappa$.
In both cases, $\Theta\in J_\kappa$, which completes the proof.
\end{proof}

\begin{lemma}\label{lem:central-part}
It holds that $\Ocal^\C(M)R_H\subseteq\cL^\C$ for every $H\in\frakz$.
\end{lemma}

\begin{proof}
Let $\phi\in\Phi$, $f\in\Ocal^\C(M)$, $X\in\fraks$, and $H\in\frakz$. By~\eqref{eq:bracket-degree} and $[R_H,R_X]=0$,
$$
[\phi R_H,fR_X]=\phi(R_Hf)R_X-f(R_X\phi)R_H.
$$
The left-hand side belongs to $\cL^\C$, since $\phi R_H\in\sX_1^\C$ and $fR_X\in\cL^\C$ by Lemma~\ref{lem:semisimple-part}, and so does the first term on the right-hand side, by Lemma~\ref{lem:semisimple-part}. Thus, $f(R_X\phi)R_H\in\cL^\C$.
Since $\fraks\cdot\Phi=\Phi_S$ by~\eqref{eq:s-Phi}, it follows that $f\phi R_H\in\cL^\C$ for all $f\in\Ocal^\C(M)$ and $\phi\in\Phi_S$.
Now, fix a $K$-invariant inner product on $\mathbb V$. Then, $\mathbb V_0$ is the orthogonal complement of $\mathbb V_S$, since the orthogonal complement of $\mathbb V_S$ is a $K$-submodule isomorphic to $\mathbb V/\mathbb V_S\approx\mathbb V_0$, on which $S$ acts trivially, and $\mathbb V=\mathbb V_S\oplus\mathbb V_0$.
Let $b_1,\ldots,b_q$ be an orthonormal basis of $\mathbb V_S$, and $\phi_j\in\Phi_S$ be the restriction to $M$ of $y\mapsto\langle y,b_j\rangle$.
Then, $\sum_{j=1}^q\phi_j^2$ is the restriction to $M$ of the $K$-invariant function $y\mapsto|y_S|^2$, where $y_S$ is the orthogonal projection of $y$ onto $\mathbb V_S$; since $K$ acts transitively on $M$, it is a constant $c_0\geq0$.
If $c_0=0$, then $M\subseteq \mathbb V_0$, and $S$ acts trivially on $M$, which contradicts the standing assumption.
Thus, $c_0>0$, and $fR_H=c_0^{-1}\sum_{j=1}^q(f\phi_j)\phi_jR_H\in\cL^\C$.
\end{proof}

With the two lemmas above, we now prove Theorem~\ref{thm:reductive}:

\begin{proof}[Proof of Theorem~\ref{thm:reductive}]
By Lemmas~\ref{lem:semisimple-part} and~\ref{lem:central-part}, $\cL^\C\supseteq\Ocal^\C(M)\sX_0^\C$, which is $\xalg^\C(M)$ by Lemma~\ref{lem:tangent}.
Since $\cL\subseteq\xalg(M)$, taking the $\sigma$-fixed points gives $\cL=\xalg(M)$.
Since $\cL\subseteq\operatorname{Lie}(\sX_0+\sX_1)\subseteq\xalg(M)$, the last assertion follows.
\end{proof}

\subsection{The case where the semi-simple part acts trivially}\label{ssec:tori}

 Throughout this subsection, we assume that $K$ is connected and that $S$ acts trivially on $M$.
We first show that $M$ is a torus on which $Z$ acts transitively:

\begin{lemma}\label{lem:torus-structure}
 Suppose that $K$ is connected and that $S$ acts trivially on $M$.
Then, the torus $Z$ acts transitively on $M$, and $zZ_x\mapsto z\cdot x$ is a diffeomorphism from the torus $Z/Z_x$ onto $M$, for any $x\in M$.
\end{lemma}

\begin{proof}
 Note that $K\cdot x=ZS\cdot x=Z\cdot x$.
The stabilizer $Z_x$ is a closed subgroup of the torus $Z$, so $Z/Z_x$ is a compact, connected, abelian Lie group, i.e., a torus, and $zZ_x\mapsto z\cdot x$ is a diffeomorphism from $Z/Z_x$ onto $M$.
\end{proof}

We fix $x\in M$ and identify $M$ with $Z/Z_x$ as in Lemma~\ref{lem:torus-structure}, and $Z/Z_x$ with $\T^m:=\R^m/2\pi\Z^m$ via an isomorphism of Lie groups, where $m=\dim M$.
Let $\theta=(\theta_1,\ldots,\theta_m)$ be the resulting angle coordinates and $\partial_j:=\partial/\partial\theta_j$.
Then, $K$ acts on $\T^m$ by translations, $\sX_0=\operatorname{span}\{\partial_1,\ldots,\partial_m\}$,
and $\Ocal(M)$ is the algebra of real trigonometric polynomials in $\theta$, since the $K$-finite functions are the finite sums of characters.
For $\chi\in\Z^m$ and $u\in\C^m$, we let
$$
    e_\chi:=e^{\mathrm{i}\chi\cdot\theta}
    \quad\mbox{and}\quad
    u\cdot\partial:=\sum_{j=1}^mu_j\partial_j,
$$
where $\chi\cdot\theta:=\sum_j\chi_j\theta_j$ and, more generally, $u\cdot u':=\sum_ju_ju'_j$ for $u,u'\in\C^m$.
By Lemma~\ref{lem:tangent}, $\xalg(M)$ is the Lie algebra $\xtrig(\T^m)$ of {\it trigonometric} vector fields, i.e., of vector fields $\sum_jf_j\partial_j$ with real trigonometric polynomials $f_j$,
and $\xtrig(\T^m)^\C$ is spanned by the vector fields $e_\chi u\cdot\partial$, for which
\begin{equation}\label{eq:torus-bracket}
[e_\chi u\cdot\partial,e_{\chi'}u'\cdot\partial]=\mathrm{i}\,e_{\chi+\chi'}\bigl((u\cdot\chi')u'-(u'\cdot\chi)u\bigr)\cdot\partial.
\end{equation}

For a finite set $C\subseteq\Z^m\setminus\{0\}$, let $\mathfrak l_C$ be the complex Lie algebra generated by the vector fields $\partial_j$ and $e_\chi\partial_j$, for $\chi\in C$ and $j=1,\ldots,m$.
The space $\Phi:=\mathcal A_1^\C$ is a finite-dimensional space of trigonometric polynomials that is invariant under translations, and hence, it is spanned by the characters $e_\chi$ that it contains.
Let $C\subseteq\Z^m\setminus\{0\}$ be the set of nonzero $\chi$ with $e_\chi\in\Phi$.
Then, $C=-C$, since $\Phi$ is closed under complex conjugation, and $C$ generates the group $\Z^m$,
since $\Phi$ generates the algebra of all trigonometric polynomials by Lemma~\ref{lem:kfinite}.
Since $(\sX_0+\sX_1)^\C$ is spanned by the vector fields $\partial_j$ and $e_\chi\partial_j$, for $\chi\in C$ and $j=1,\ldots,m$, we have that $\operatorname{Lie}(\sX_0+\sX_1)^\C=\mathfrak l_C$, and the following result covers the case where $S$ acts trivially:

\begin{theorem}\label{thm:torus}
Let $m\in\N$, and $C\subseteq\Z^m\setminus\{0\}$ be a finite set with $C=-C$ that generates the group $\Z^m$.
If $(m,C)\neq(1,\{\pm1\})$, then $\mathfrak l_C=\xtrig(\T^m)^\C$.
If $(m,C)=(1,\{\pm1\})$, then $\mathfrak l_C=\operatorname{span}\{\partial_1,e^{\mathrm{i}\theta_1}\partial_1,e^{-\mathrm{i}\theta_1}\partial_1\}$, which is isomorphic to $\mathfrak{sl}_2(\C)$.
\end{theorem}

The proof uses the interpolation argument in the proof of Kuranishi's theorem~\cite{Kuranishi49}, \cite[Theorem~6]{Kuranishi51}, which will be used repeatedly in Sections~\ref{sec:two-generators} and~\ref{sec:spheres} as well.
We state it for an arbitrary complex Lie algebra:

\begin{lemma}\label{lem:interpolation}
Let $L$ be a complex Lie algebra, $W_1\in L$, and $W_2=\sum_{j=1}^nV_j$ with $V_j\in L$ such that $[W_1,V_j]=c_jV_j$ for pairwise distinct complex numbers $c_1,\ldots,c_n$. Then, every $V_j$ belongs to the complex Lie algebra generated by $W_1$ and $W_2$.
\end{lemma}

\begin{proof}
We have that $\operatorname{ad}_{W_1}^pW_2=\sum_{j=1}^nc_j^pV_j$ for every $p\in\N_0$.
For $j=1,\ldots,n$, let $q_j(z):=\prod_{k\neq j}(z-c_k)/(c_j-c_k)$. Then, $q_j(\operatorname{ad}_{W_1})W_2=V_j$, and the left-hand side is a linear combination of $W_2$ and its iterated Lie brackets with $W_1$.
\end{proof}

With the lemma above, we now prove Theorem~\ref{thm:torus}:

\begin{proof}[Proof of Theorem~\ref{thm:torus}]
For $\chi\in\Z^m$, let $L_\chi:=\{u\in\C^m\mid e_\chi u\cdot\partial\in\mathfrak l_C\}$.
By~\eqref{eq:torus-bracket}, $[t\cdot\partial,e_\chi u\cdot\partial]=\mathrm{i}\,(\chi\cdot t)\,e_\chi u\cdot\partial$ for $t\in\R^m$.
Let $A=\sum_\chi e_\chi u_\chi\cdot\partial$ be an element of $\mathfrak l_C$, where the sum is over a finite set of frequencies,
and choose $t\in\R^m$ such that the numbers $\chi\cdot t$ are pairwise distinct for these frequencies.
Since $t\cdot\partial\in\mathfrak l_C$, it follows from Lemma~\ref{lem:interpolation}, applied to $W_1=t\cdot\partial$ and $W_2=A$, that each $e_\chi u_\chi\cdot\partial$ belongs to $\mathfrak l_C$. Thus,
$$
  \mathfrak l_C=\bigoplus_{\chi\in\Z^m}e_\chi L_\chi\cdot\partial.
$$
Let $\Gamma:=\{\chi\in\Z^m\mid L_\chi=\C^m\}$. By hypothesis, $C\cup\{0\}\subseteq\Gamma$.

First, $\Gamma=-\Gamma$. To wit, the diffeomorphism $\theta\mapsto-\theta$ of $\T^m$ induces an automorphism $\iota$ of the Lie algebra of complex vector fields with $\iota(e_\chi u\cdot\partial)=-e_{-\chi}u\cdot\partial$. Since $C=-C$,
    the automorphism $\iota$ preserves the set of generators of $\mathfrak l_C$ up to sign, so $\iota(\mathfrak l_C)=\mathfrak l_C$ and $L_{-\chi}=L_\chi$ for all $\chi$.

Second, if $\chi,\chi'\in\Gamma\setminus\{0\}$ are distinct, then $\chi+\chi'\in\Gamma$.
If $\chi+\chi'=0$, there is nothing to show. Otherwise, by~\eqref{eq:torus-bracket}, $L_{\chi+\chi'}$ contains the span $B$ of the vectors $(u\cdot\chi')u'-(u'\cdot\chi)u$, for $u,u'\in\C^m$.
If $\chi$ and $\chi'$ are linearly independent, we choose $u$ with $u\cdot\chi'=1$ and $u\cdot\chi=0$, and see that $B$ contains every $u'$ with $u'\cdot\chi=0$. Similarly, $B$ contains every $u$ with $u\cdot\chi'=0$. These two hyperplanes are distinct, so $B=\C^m$.
If $\chi'=c\chi$ for some $c\in\R$, then $c\neq1$, and taking $u'=u$ with $u\cdot\chi=1$, we obtain $(c-1)u\in B$; such vectors $u$ span $\C^m$, so $B=\C^m$.

Third, if for every $\chi\in\Gamma\setminus\{0\}$, there exists $\chi'\in\Gamma\setminus\{0,\chi,-\chi\}$, then $\Gamma$ is closed under addition.
To wit, by the first two claims, it suffices to show that $2\chi\in\Gamma$ for $\chi\in\Gamma\setminus\{0\}$. With $\chi'$ as above,
$\chi+\chi'$ and $\chi-\chi'$ are distinct nonzero elements of $\Gamma$, by the second claim applied to $(\chi,\chi')$ and $(\chi,-\chi')$, and hence, $2\chi=(\chi+\chi')+(\chi-\chi')\in\Gamma$.
Under the hypothesis of the third claim, since $\Gamma$ contains $C=-C$ and is closed under addition, it contains the group generated by $C$, which is $\Z^m$, and hence, $\mathfrak l_C=\xtrig(\T^m)^\C$.

It remains to verify the hypothesis of the third claim unless $m=1$ and $C=\{\pm1\}$.
If $m\geq2$, then for every $\chi\neq0$, some $\chi'\in C$ is not proportional to $\chi$, since $C$ generates $\Z^m$, and $\chi'\in\Gamma\setminus\{0,\pm\chi\}$.
If $m=1$ and $C\neq\{\pm1\}$, then $C$ contains two elements $0<a<b$, since $C=-C$ generates $\Z$, and for every $\chi\neq0$, one of $a$ and $b$ belongs to $\Gamma\setminus\{0,\pm\chi\}$.
Finally, if $m=1$ and $C=\{\pm1\}$, then $[e^{\mathrm{i}\theta_1}\partial_1,e^{-\mathrm{i}\theta_1}\partial_1]=-2\mathrm{i}\,\partial_1$ and $[\partial_1,e^{\pm\mathrm{i}\theta_1}\partial_1]=\pm\mathrm{i}\,e^{\pm\mathrm{i}\theta_1}\partial_1$. Thus,
  the span of $\partial_1$ and $e^{\pm\mathrm{i}\theta_1}\partial_1$ is closed under the Lie bracket, and it is isomorphic to $\mathfrak{sl}_2(\C)$.
\end{proof}

\subsection{Proof of Theorem~\ref{thm:every-embedding}}\label{ssec:proof-every-embedding}

With the results above, we now establish Theorem~\ref{thm:every-embedding}:

\begin{proof}[Proof of Theorem~\ref{thm:every-embedding}]
By Lemma~\ref{lem:subgroup}, we can assume that $K$ is connected.
If $S$ acts nontrivially on $M$, then $\xalg(M)=\operatorname{Lie}(\sX_0+\sX_1)$ by Theorem~\ref{thm:reductive}.
Moreover, the exceptional case does not occur: in that case, $S$ would act on the circle $M$ through the abelian group of rotations,
    and since $S=[S,S]$, it would act trivially.
If $S$ acts trivially on $M$, then, with the notation of Subsection~\ref{ssec:tori}, $\operatorname{Lie}(\sX_0+\sX_1)^\C=\mathfrak l_C$.
If $(m,C)\neq (1,\{\pm1\})$, then $\mathfrak l_C=\xalg^\C (M)$ by Theorem~\ref{thm:torus}, and taking the $\sigma$-fixed points gives $\xalg(M)=\operatorname{Lie}(\sX_0+\sX_1)$.
In this case, the exceptional case does not occur either.
To wit, suppose that $m=1$, and let $\theta'$ be an angle coordinate on $M$ in which $K$ acts by rotations.
Then, $K$ acts transitively by translations of $\R/2\pi\Z$ in both coordinates $\theta_1$ and $\theta'$, so the change of coordinates is, up to a translation, a continuous automorphism of $\R/2\pi\Z$; thus, $\theta'=\pm\theta_1+c$ for some $c\in\R$.
Hence, if $\mathcal A_1\subseteq\operatorname{span}\{1,\cos\theta',\sin\theta'\}$, then every function in $\mathcal A_1$ is a linear combination of $1$ and $e^{\pm\mathrm{i}\theta_1}$, so $C\subseteq\{\pm1\}$, and hence, $C=\{\pm1\}$, since $C$ generates $\Z$, which contradicts $(m,C)\neq(1,\{\pm1\})$.
If $(m,C)=(1,\{\pm1\})$, then $M$ is a circle, $\mathcal A_1\subseteq\operatorname{span}\{1,\cos\theta_1,\sin\theta_1\}$, and by Theorem~\ref{thm:torus}, $\operatorname{Lie}(\sX_0+\sX_1)$ is the real span of $\partial_1$, $\cos\theta_1\,\partial_1$, and $\sin\theta_1\,\partial_1$, which is isomorphic to $\mathfrak{sl}_2(\R)$.
\end{proof}

\section{Two generators and the proof of Theorem~\ref{thm:main}}\label{sec:two-generators}

In this section, we prove Theorem~\ref{thm:main}.
By Lemma~\ref{lem:subgroup}, we can assume that $K$ is connected, and we use the notation of Section~\ref{sec:degree-growth}.
If the semi-simple part $S$ of $K$ acts trivially on $M$, then $M$ is a torus on which $K$ acts by translations (see Subsection~\ref{ssec:tori}),
and two generators are given in Proposition~\ref{prop:torus-two} below.
We thus assume, until Subsection~\ref{ssec:torus-two}, that $S$ acts nontrivially on $M$.

The two generators will be $W_1=R_T$ and $W_2=R_E+w$, where $T$ belongs to a maximal abelian subalgebra of $\frakk$,
  $E\in\fraks$ is the sum of nonzero simple root vectors and their complex conjugates, and $w\in\sX_1$.
By Lemma~\ref{lem:interpolation}, $\operatorname{Lie}\{W_1,W_2\}^\C$ contains the weight components of $R_E$ and of $w$, provided that the eigenvalues of $\operatorname{ad}_{R_T}$ on these components are pairwise distinct.
In this case, it also contains the vector fields $R_X$, for $X\in\fraks^\C$, since the weight components of $R_E$ are the vector fields $R_X$,
    where $X$ is a nonzero simple root vector or its complex conjugate.
By Theorem~\ref{thm:reductive}, it then suffices that the weight components of $w$ generate $\sX_1^\C$ as an $\fraks^\C$-module.
However, an irreducible summand may occur in $\sX_1^\C$ with a multiplicity larger than the number of the weights at which the components of $w$ can be separated from those of $R_E$; then,
    the weight components of $w$ do not generate $\sX_1^\C$ for any $w$, and the above argument does not apply (see Example~\ref{rem:W-fails}).
We show that this does not happen for suitable equivariant embeddings, which we construct in Subsections~\ref{ssec:weights} and~\ref{ssec:spread}; the generators are given in Theorem~\ref{thm:two-generators} in Subsection~\ref{ssec:separation}, and tori are treated in Subsection~\ref{ssec:torus-two}.

\subsection{Weights and multiplicities}\label{ssec:weights}

We first fix notation for weights and bound the multiplicities of the irreducible summands of $\sX_1^\C$.
Let $\frakt$ be a maximal abelian subalgebra of $\frakk$, i.e., the Lie algebra of a maximal torus of $K$, and $r:=\dim\frakt$.
Then, $\frakt=\frakz\oplus\frakh$, where $\frakh:=\frakt\cap\fraks$ is a maximal abelian subalgebra of $\fraks$.
On every finite-dimensional complex $K$-module $U$, the elements of $\frakt$ act by commuting diagonalizable operators with purely imaginary eigenvalues.
For $\nu\in(\frakt^\C)^*$, we let $U[\nu]:=\{u\in U\mid X\cdot u=\nu(X)u\mbox{ for all }X\in\frakt\}$, and we call $\nu$ a {\it weight} of $U$ if $U[\nu]\neq0$.
Since every algebraic vector field is contained in a finite-dimensional $K$-invariant subspace of $\xalg(M)$,
  the above applies to $\xalg^\C(M)$, on which $X\in\frakk$ acts by $V\mapsto[R_X,V]$ by~\eqref{eq:R-derivatives}.
Let $\Delta\subseteq(\frakt^\C)^*$ be the set of roots, i.e., of nonzero weights of the adjoint module $\frakk^\C$; the roots vanish on $\frakz$, and the root spaces $\frakk^\C[\alpha]$ are contained in $\fraks^\C$.
We fix a positive system $\Delta^+\subseteq\Delta$, and let $\Pi$ be the set of simple roots and $\frakn^+:=\bigoplus_{\alpha\in\Delta^+}\frakk^\C[\alpha]$.
The complex conjugation $\sigma$ of $\frakk^\C$ maps $\frakk^\C[\alpha]$ onto $\frakk^\C[-\alpha]$, and for a real $K$-module $U$, it maps $U^\C[\nu]$ onto $U^\C[-\nu]$, since the weights are purely imaginary on $\frakt$.

Every finite-dimensional complex $K$-module is completely reducible, and the irreducible ones are classified by their highest weights with respect to $\Delta^+$ (see, e.g.,~\cite[Chapter~VI]{BtD}).
For a dominant integral element $\lambda\in (\frakt^\C)^*$, we let $V_\lambda$ be the irreducible $K$-module with $\lambda$ the highest weight, and for a finite-dimensional complex $K$-module $U$, we let $c_\lambda(U)$ be the multiplicity of $V_\lambda$ in $U$.
Since $V_\lambda$ contains, up to scalars, a unique nonzero vector annihilated by $\frakn^+$, which has weight $\lambda$, we have that $c_\lambda(U)=\dim U^{\frakn^+}[\lambda]$, where $U^{\frakn^+}$ is the space of vectors in $U$ annihilated by $\frakn^+$.
We fix an $\operatorname{Ad}(K)$-invariant inner product on $\frakk$; the subspaces $\frakz$ and $\fraks$, and hence $\frakz$ and $\frakh$, are orthogonal, since $\langle Z,[X,Y]\rangle=-\langle[X,Z],Y\rangle=0$ for $Z\in\frakz$ and $X,Y\in\frakk$.
It induces an inner product $\langle\cdot,\cdot\rangle$, with norm $\|\cdot\|$, on the real span of the weights, which are purely imaginary on $\frakt$.
For $\nu\in(\frakt^\C)^*$, we let $\nu_S\in(\frakt^\C)^*$ be given by $\nu_S|_{\frakh}=\nu|_{\frakh}$ and $\nu_S|_{\frakz}=0$, so that $\|\nu_S\|\leq\|\nu\|$ and $\alpha_S=\alpha$ for $\alpha\in\Delta$.
Let $\rho_{\max}:=\max_{\alpha\in\Delta}\|\alpha\|$.
We need the following bound on multiplicities in tensor products:

\begin{lemma}\label{lem:multiplicity}
Let $U$ be a finite-dimensional complex $K$-module, and $\mu,\lambda\in(\frakt^\C)^*$ be highest weights.  Then,
$$
c_\lambda(V_\mu\otimes U)\leq\dim U[\lambda-\mu].
$$
\end{lemma}

\begin{proof}
Let $V_\mu^*$ be the dual module, and $v^*\in V_\mu^*$ be a nonzero vector of weight $-\mu$, which is a lowest weight vector of $V_\mu^*$, i.e., it is annihilated by $\frakk^\C[-\alpha]$ for all $\alpha\in\Delta^+$.
By the highest weight theory, applied with respect to the positive system $-\Delta^+$, the module $V_\mu^*$ is spanned by the vectors $Y_1\cdots Y_jv^*$, for $j\in\N_0$ and $Y_1,\ldots,Y_j \in \frakn^+$ (see, e.g.,~\cite[Section~20.2]{HumphreysLie}).
We identify $V_\mu\otimes U$ with the $K$-module of linear maps from $V_\mu^*$ to $U$ via $(v\otimes u)(f):=f(v)u$, for $v\in V_\mu$, $u\in U$, and $f\in V_\mu^*$; 
then, $X\in\frakk$ acts on such a linear map $A$ by $(X\cdot A)(f)=X\cdot A(f)-A(X\cdot f)$, for $f\in V_\mu^*$, since $(X\cdot f)(v)=-f(X\cdot v)$.
If $A$ is annihilated by $\frakn^+$, then $A(Y\cdot f)=Y\cdot A(f)$ for all $Y\in\frakn^+$ and $f\in V_\mu^*$, and hence, $A$ is determined by $A(v^*)$.
If, moreover, $A$ has weight $\lambda$, then $A(v^*)\in U[\lambda-\mu]$, since $H\cdot A(v^*)=(H\cdot A)(v^*)+A(H\cdot v^*)=(\lambda-\mu)(H)A(v^*)$ for $H\in\frakt$.
Thus, the linear map $A\mapsto A(v^*)$ from $(V_\mu\otimes U)^{\frakn^+}[\lambda]$ to $U[\lambda-\mu]$ is injective.
\end{proof}

We apply Lemma~\ref{lem:multiplicity} to the space $\sX_1^\C$.
For an equivariant embedding $(\mathbb V,v)$ of $M$, let $\Lambda(\mathbb V)$ be the finite multiset of highest weights such that
$$
  (\mathbb V^*)^\C\approx\bigoplus_{\mu\in\Lambda(\mathbb V)}V_\mu.
$$
Since $k\cdot((\ell|_M)R_X)=((k\cdot\ell)|_M)R_{\operatorname{Ad}_kX}$, the map $\ell\otimes X\mapsto(\ell|_M)R_X$ from $(\mathbb V^*)^\C\otimes\frakk^\C$ to $\sX_1^\C$ is a surjective homomorphism of $K$-modules.
Since $\dim\frakk^\C[0]=r$, $\dim\frakk^\C[\alpha]=1$ for $\alpha\in\Delta$, and $\frakk^\C[\beta]=0$ otherwise, Lemma~\ref{lem:multiplicity} yields that for every highest weight $\lambda$,
\begin{equation}\label{eq:multiplicity-bound}
c_\lambda(\sX_1^\C)\leq\sum_{\mu\in\Lambda(\mathbb V)}\dim\frakk^\C[\lambda-\mu]
\leq r\cdot\#\{\mu\in\Lambda(\mathbb V)\mid\lambda-\mu\in\Delta\cup\{0\}\}.
\end{equation}

\subsection{Spread-out embeddings}\label{ssec:spread}

We now construct equivariant embeddings for which the right-hand side of~\eqref{eq:multiplicity-bound} is small and the highest weights in $\Lambda(\mathbb V)$ are large.
For an irreducible complex $K$-module $U$ with highest weight $\mu$, let $\mu^*$ be the highest weight of the dual module $U^*$, and $\ell_0\in U^*$ be a highest weight vector.
For $k\in\N$, let $N_k(U)$ be the span of the homogeneous polynomial functions $y\mapsto(h\cdot\ell_0)(y)^k$ on $U$, for $h\in K$.
It is the $K$-submodule of the space of polynomial functions on $U$ generated by the highest weight vector $\ell_0^k$, so it is irreducible with highest weight $k\mu^*$,
  and its dual module $N_k(U)^*$ is irreducible with highest weight $k\mu$.
Let $\mathrm{ev}_k\colon U\to N_k(U)^*$ be given by $\mathrm{ev}_k(y)(f):=f(y)$ for $f\in N_k(U)$.
This map is $K$-equivariant, since $\mathrm{ev}_k(h\cdot y)(f)=f(h\cdot y)=(h^{-1}\cdot f)(y)=(h\cdot\mathrm{ev}_k(y))(f)$.
We have the following result:

\begin{lemma}\label{lem:spread}
Suppose that $S$ acts nontrivially on $M$, and let $R>0$.
Then, there is an equivariant embedding $(\mathbb V',v')$ of $M$ such that
\begin{enumerate}
\item $\|\mu_S\|\geq R$ for every $\mu\in\Lambda(\mathbb V')$;
\item for every $\lambda\in(\frakt^\C)^*$, at most two elements $\mu$ of $\Lambda(\mathbb V')$, counted with multiplicity, satisfy $\lambda-\mu\in\Delta\cup\{0\}$.
\end{enumerate}
\end{lemma}

\begin{proof}
Let $(\mathbb V,v)$ be an equivariant embedding of $M$. We proceed in three steps.

\xc{Step~1: removal of $S$-invariant vectors.}
As in Subsection~\ref{ssec:reductive}, we write $\mathbb V=\mathbb V_S\oplus \mathbb V_0$, and $v=v_S+v_0$ accordingly.
Then, $v_S\neq0$; otherwise, $v\in \mathbb V_0$, and since $S$ is normal, $s\cdot(g\cdot v)=g\cdot((g^{-1}sg)\cdot v)=g\cdot v$ for all $s\in S$ and $g\in K$, i.e., $S$ acts trivially on $M$.
 We choose a decomposition of $\mathbb V_S$ into irreducible $K$-submodules and a summand $\mathbb U$ of this decomposition such that the component $y$ of $v_S$ in $\mathbb U$ is nonzero; since the projection of $\mathbb V_S$ onto $\mathbb U$ along the other summands is $K$-equivariant, every element of $K$ that fixes $v_S$ also fixes $y$.
Let $\mathbb V'':= \mathbb V_S\oplus(\mathbb U\otimes \mathbb V_0)$, with the diagonal action on the tensor product, and $v'':=v_S+y\otimes v_0$.
If $g\in K_v$, then $g$ fixes $v_S$ and $v_0$, and hence, $y$; thus, $g\in K_{v''}$.
Conversely, if $g\in K_{v''}$, then $g$ fixes $v_S$, and hence, $y$, and then $y\otimes(g\cdot v_0)=y\otimes v_0$ implies that $g\cdot v_0=v_0$, since $y\neq0$.
Thus, $K_{v''}=K_v$, and $(\mathbb V'',v'')$ is an equivariant embedding of $M$.
 It does not contain any nonzero $S$-invariant vector.
To wit, $\mathbb V_S$ does not contain any nonzero $S$-invariant vector, since $\mathbb V_S\cap\mathbb V_0=0$, and hence, neither does its submodule $\mathbb U$; since $S$ acts trivially on $\mathbb V_0$, the $S$-module $\mathbb U\otimes\mathbb V_0$ is isomorphic to a direct sum of copies of $\mathbb U$, so it does not contain any nonzero $S$-invariant vector either.

\xc{Step~2: Cartan powers.}
We decompose $\mathbb V''=\mathbb V_1\oplus\cdots\oplus\mathbb V_p$ into irreducible $K$-submodules, and write $v''=v_1+\cdots+v_p$ accordingly; removing from $\mathbb V''$ the summands with $v_j=0$, which does not change $K_{v''}$, we can assume that $v_j\neq0$ for all $j$.
Then, $K_{v''}=\bigcap_jK_{v_j}$, and $S$ acts nontrivially on every $\mathbb V_j$, since $\mathbb V''$ does not contain any nonzero $S$-invariant vector by Step~1.
The modules $N_k(U)$ have been defined for irreducible complex $K$-modules $U$, whereas the $K$-modules $\mathbb V_j$ are real.
We therefore choose, for each $j$, an irreducible complex $K$-submodule $U_j$ of $\mathbb V_j^\C$, and we let $\mu_j$ be its highest weight and $\pi_j\colon\mathbb V_j\to U_j$ be the restriction to $\mathbb V_j$ of a $K$-equivariant projection of $\mathbb V_j^\C$ onto $U_j$.
The real-linear map $\pi_j$ is $K$-equivariant and nonzero, since $\mathbb V_j$ spans $\mathbb V_j^\C$ over $\C$, and it is injective, since $\mathbb V_j$ is irreducible.
Thus, $x_j:=\pi_j(v_j)\neq0$ and $K_{x_j}=K_{v_j}$.
Since $U_j+\sigma(U_j)$ is the complexification of a nonzero $K$-submodule of $\mathbb V_j$, we have that $\mathbb V_j^\C=U_j+\sigma(U_j)$; since $S$ acts nontrivially on $\mathbb V_j$ and $\sigma$ commutes with the action of $K$, it follows that $S$ acts nontrivially on $U_j$.
Hence, $(\mu_j)_S\neq0$: the center $\frakz$ acts on the irreducible $K$-module $U_j$ by scalars, so $U_j$ is an irreducible $\fraks^\C$-module with highest weight $\mu_j|_{\frakh^\C}$, and if $\mu_j|_{\frakh^\C}=0$, then $U_j$ is the trivial $\fraks^\C$-module, on which $S$ acts trivially.

Let $k,k'\in\N$ be coprime, and let $g\in K$ satisfy $\mathrm{ev}_k(g\cdot x_j)=\mathrm{ev}_k(x_j)$ and $\mathrm{ev}_{k'}(g\cdot x_j)=\mathrm{ev}_{k'}(x_j)$,
    where $\mathrm{ev}_k$ and $\mathrm{ev}_{k'}$ are the maps defined above for $U=U_j$. We claim that $g\cdot x_j=x_j$.
To wit, let $h\in K$, $a:=(h\cdot\ell_0)(x_j)$, and $b:=(h\cdot\ell_0)(g\cdot x_j)$. Then, $a^k=b^k$ and $a^{k'}=b^{k'}$. If $a\neq0$, then $b/a$ is both a $k$th and a $k'$th root of unity, and hence, $b=a$; if $a=0$, then $b=0$.
Thus, $(h\cdot\ell_0)(g\cdot x_j-x_j)=0$ for all $h\in K$, and since the vectors $h\cdot\ell_0$ span $U_j^*$, we have that $g\cdot x_j=x_j$.
Conversely, every element of $K_{x_j}$ fixes $\mathrm{ev}_k(x_j)$ and $\mathrm{ev}_{k'}(x_j)$ by equivariance.

Now, given $k_j,k_j'\in\N$ with $\gcd(k_j,k_j')=1$, for $j=1,\ldots,p$, let
$$
  \mathbb V':=\bigoplus_{j=1}^p\bigl(N_{k_j}(U_j)^*\oplus N_{k_j'}(U_j)^*\bigr)_\R
  \quad\mbox{and}\quad
  v':=\sum_{j=1}^p\bigl(\mathrm{ev}_{k_j}(x_j)+\mathrm{ev}_{k_j'}(x_j)\bigr),
$$
 where, for a complex $K$-module $U'$, we let $U'_\R$ be its realification, i.e., $U'$ regarded as a real $K$-module; thus, $\mathbb V'$ is a real $K$-module, as required for an equivariant embedding.
By the above claim, $K_{v'}=\bigcap_jK_{x_j}=\bigcap_jK_{v_j}=K_{v''}=K_v$, so $(\mathbb V',v')$ is an equivariant embedding of $M$.

\xc{Step~3: the weights.}
For a complex $K$-module $U'$, we have that $(U'_\R)^\C\approx U'\oplus\overline{U'}$, and the complex conjugate module $\overline{U'}$ is isomorphic to $U'^*$ via a $K$-invariant Hermitian inner product. Also, $\mathbb V'^*\approx \mathbb V'$ via a $K$-invariant inner product. Thus,
$$
\Lambda(\mathbb V')=\{k_j\mu_j,\ k_j\mu_j^*,\ k_j'\mu_j,\ k_j'\mu_j^*\mid j=1,\ldots,p\},
$$
as a multiset.
For an irreducible complex $K$-module $U$ with highest weight $\mu$, the restriction of $\mu^*$ to $\frakh^\C$ is $-w_{\mathrm{long}}(\mu|_{\frakh^\C})$, where $w_{\mathrm{long}}$ is the longest element of the Weyl group of $\fraks^\C$ with respect to $\Pi$, i.e., the element that maps $\Delta^+$ onto $-\Delta^+$ (see, e.g.,~\cite[Exercise~21.6]{HumphreysLie}).
Since the Weyl group acts isometrically on the real span of the weights, it follows that $\|(\mu^*)_S\|=\|\mu_S\|$.
Thus, $\|(k\mu_j^*)_S\|=\|(k\mu_j)_S\|=k\tau_j$, where $\tau_j:=\|(\mu_j)_S\|>0$.
We fix $m_0\in\N$ with $m_0\min_j\tau_j>2\rho_{\max}$, and choose primes $k_1,\ldots,k_p>m_0$ such that, for all $j$,
$$
    k_j\tau_j\geq R
    \quad\mbox{and}\quad
    k_j\tau_j>(k_{j'}+m_0)\tau_{j'}+2\rho_{\max}\ \mbox{ for all } j'<j,
$$
and let $k_j':=k_j+m_0$. Then, $\gcd(k_j,k_j')=\gcd(k_j,m_0)=1$, and the first condition of the lemma holds.
For the second condition, note that the $2p$ numbers $k\tau_j$, for $j=1,\ldots,p$ and $k\in\{k_j,k_j'\}$, are pairwise at distance greater than $2\rho_{\max}$.
Now, let $\nu,\nu'\in\Lambda(\mathbb V')$ be such that $\lambda-\nu$ and $\lambda-\nu'$ belong to $\Delta\cup\{0\}$. Then,
    $\|\nu_S-\nu'_S\|\leq\|(\lambda-\nu)_S\|+\|(\lambda-\nu')_S\|\leq2\rho_{\max}$.
If $\nu\in\{k\mu_j,k\mu_j^*\}$ and $\nu'\in\{k''\mu_{j'},k''\mu_{j'}^*\}$, then $|k\tau_j-k''\tau_{j'}|\leq\|\nu_S-\nu'_S\|\leq2\rho_{\max}$, and hence, $(j,k)=(j',k'')$.
Thus, all such $\nu$ belong to the pair $\{k\mu_j,k\mu_j^*\}$ for a single $(j,k)$, which proves the second condition.
\end{proof}

\subsection{Separation of components}\label{ssec:separation}

We now formulate the condition on multiplicities described at the beginning of this section.

\begin{definition}\label{def:W}
For a highest weight $\lambda$, let
$$
a(\lambda):=\sum_{\nu\notin\Pi\cup(-\Pi)}\dim V_\lambda[\nu].
$$
An equivariant embedding $(\mathbb V,v)$ of $M$ is said to satisfy {\it condition~(W)} if
\begin{equation}\label{eq:condition-W}
c_\lambda(\sX_1^\C)\leq a(\lambda)\quad\mbox{for every highest weight }\lambda.
\end{equation}
\end{definition} 
Condition~(W) states that, for every highest weight $\lambda$, the multiplicity of $V_\lambda$ in $\sX_1^\C$ does not exceed the number of weights of $V_\lambda$, counted with multiplicity, that are different from the weights $\pm\alpha$, for $\alpha\in\Pi$, of the weight components of $R_E$.
We have the following result:

\begin{lemma}\label{lem:W}
There is $R_0>0$, depending only on $\frakk$, such that every equivariant embedding $(\mathbb V',v')$ as in Lemma~\ref{lem:spread} with $R\geq R_0$ satisfies condition~(W).
\end{lemma}

\begin{proof}
Let $\varpi_\alpha$, for $\alpha\in\Pi$, be the fundamental weights of $\fraks^\C$ with respect to $\frakh^\C$ and $\Pi$,
extended by zero on $\frakz$, and let $C_0:=\sum_{\alpha\in\Pi}\|\varpi_\alpha\|$.
We show that $R_0:=\rho_{\max}+4rC_0$ has the required property.
Let $\lambda$ be a highest weight with $c_\lambda(\sX_1^\C)\neq0$.
By~\eqref{eq:multiplicity-bound} and the second condition of Lemma~\ref{lem:spread}, $c_\lambda(\sX_1^\C)\leq2r$, and there is $\mu\in\Lambda(\mathbb V')$ with $\lambda-\mu\in\Delta\cup\{0\}$, so that $\|\lambda_S\|\geq\|\mu_S\|-\rho_{\max}\geq R-\rho_{\max}\geq4rC_0$.
For $\alpha\in\Pi$, let $n_\alpha:=2\langle\lambda_S,\alpha\rangle/\|\alpha\|^2$, which is a nonnegative integer since $\lambda$ is dominant, and let $n:=\max_{\alpha\in\Pi}n_\alpha$.
Since $\lambda_S=\sum_{\alpha\in\Pi}n_\alpha\varpi_\alpha$, we have that $\|\lambda_S\|\leq nC_0$, and hence, $n\geq4r$.
Let $\alpha\in\Pi$ be such that $n_\alpha=n$. Then, $\lambda,\lambda-\alpha,\ldots,\lambda-n\alpha$ are $n+1$ distinct weights of $V_\lambda$ (see,
    e.g.,~\cite[Section~21.3]{HumphreysLie}), and at most $2|\Pi|\leq2r$ of them belong to $\Pi\cup(-\Pi)$.
Thus, $a(\lambda)\geq n+1-2r>2r\geq c_\lambda(\sX_1^\C)$.
\end{proof}

Let $\Omega^\circ$ be the set of weights of $\sX_1^\C$ that do not belong to $\Pi\cup(-\Pi)$, and let $\sX_1^\circ$ be the space of $w\in\sX_1$ whose components in $\sX_1^\C[\nu]$, for $\nu\in\Pi\cup(-\Pi)$, vanish.
Since $\sigma$ maps $\sX_1^\C[\nu]$ onto $\sX_1^\C[-\nu]$, the complexification of $\sX_1^\circ$ is $\bigoplus_{\nu\in\Omega^\circ}\sX_1^\C[\nu]$.
Let $\Omega\subseteq(\frakt^\C)^*$ be the finite set consisting of $0$, the roots, and the weights of $\sX_1^\C$.
We fix $T\in\frakt$ such that
\begin{equation}\label{eq:separation}
\nu(T)\neq\nu'(T)\quad\mbox{whenever }\nu,\nu'\in\Omega\mbox{ are distinct};
\end{equation}
such $T$ exist, since for distinct $\nu,\nu'\in\Omega$, the condition $\nu(T)=\nu'(T)$ defines a proper hyperplane of $\frakt$, and a finite union of proper hyperplanes does not cover $\frakt$.
For each $\alpha\in\Pi$, we fix a nonzero $E_\alpha\in\frakk^\C[\alpha]$, and we let
\begin{equation}\label{eq:E}
E_+:=\sum_{\alpha\in\Pi}E_\alpha
\quad\mbox{and}\quad
E:=E_++\sigma(E_+)\in\fraks.
\end{equation}
We recall that a subset of a finite-dimensional real vector space is {\it Zariski open} if its complement is the common zero set of finitely many real polynomial functions; a nonempty Zariski open subset is dense, since a polynomial function that vanishes on a nonempty open subset is zero.
We now establish the following result:

\begin{lemma}\label{lem:separation}
Suppose that $S$ acts nontrivially on $M$ and that $(\mathbb V,v)$ satisfies condition~(W).
Then, there is a nonempty Zariski open subset $\mathcal G$ of the real vector space $\sX_1^\circ$ such that $\operatorname{Lie}\{R_T,R_E+w\}\supseteq\sX_0^{\fraks}+\sX_1$ for every $w\in\mathcal G$.
\end{lemma}

 \begin{proof}
For every $w$ in the complexification of $\sX_1^\circ$, we write $w=\sum_{\nu\in\Omega^\circ}w_\nu$ with $w_\nu\in\sX_1^\C[\nu]$, and we let
$$
\mathcal M_w:=\sum_{\nu\in\Omega^\circ}\mathcal U(\fraks^\C)\,w_\nu
$$
be the $\fraks^\C$-submodule of $\sX_1^\C$ generated by the components $w_\nu$, where the universal enveloping algebra $\mathcal U(\fraks^\C)$ acts on $\sX_1^\C$ through $X\mapsto\operatorname{ad}_{R_X}$.
We proceed in three steps.

\xc{Step~1: $\frakw^\C\supseteq\sX_0^{\fraks,\C}+\mathcal M_w$ for every $w\in\sX_1^\circ$, where $\frakw:=\operatorname{Lie}\{R_T,R_E+w\}$.}
Since $X\mapsto R_X$ is a Lie algebra homomorphism, $[R_T,R_{E_\alpha}]=\alpha(T)R_{E_\alpha}$ and $[R_T,R_{\sigma(E_\alpha)}]=-\alpha(T)R_{\sigma(E_\alpha)}$ for $\alpha\in\Pi$, and $[R_T,w_\nu]=\nu(T)w_\nu$ for $\nu\in\Omega^\circ$ by the definition of the weight spaces.
Thus,
$$
R_E+w=\sum_{\alpha\in\Pi}\bigl(R_{E_\alpha}+R_{\sigma(E_\alpha)}\bigr)+\sum_{\nu\in\Omega^\circ}w_\nu
$$
is a sum of eigenvectors of $\operatorname{ad}_{R_T}$, whose eigenvalues $\alpha(T)$ and $-\alpha(T)$, for $\alpha\in\Pi$, and $\nu(T)$, for $\nu\in\Omega^\circ$, are pairwise distinct by~\eqref{eq:separation}.
By Lemma~\ref{lem:interpolation}, every summand belongs to $\frakw^\C$; in particular, $\frakw^\C$ contains $R_{E_\alpha}$ and $R_{\sigma(E_\alpha)}$, for $\alpha\in\Pi$, and $w_\nu$, for $\nu\in\Omega^\circ$.
The root vectors $E_\alpha$ and $\sigma(E_\alpha)$, for $\alpha\in\Pi$, generate $\fraks^\C$ (see, e.g.,~\cite[Sections~8.4 and~10.2]{HumphreysLie}), so $\sX_0^{\fraks,\C}\subseteq\frakw^\C$.
Since $\frakw^\C$ is invariant under $\operatorname{ad}_{R_X}$ for $X\in\fraks^\C$, it follows that $\mathcal M_w\subseteq\frakw^\C$.

\xc{Step~2: description of $\mathcal M_w$.}
For a highest weight $\lambda$ with $c_\lambda(\sX_1^\C)\neq0$, let $I_\lambda$ be the isotypic component of type $V_\lambda$ of the $K$-module $\sX_1^\C$, and $P_\lambda\colon\sX_1^\C\to I_\lambda$ be the projection along the other isotypic components, which is a homomorphism of $K$-modules.
We fix an isomorphism of $K$-modules $I_\lambda\approx V_\lambda\otimes C_\lambda$, where $C_\lambda$ is a complex vector space of dimension $c_\lambda(\sX_1^\C)$ on which $K$ acts trivially, and for each $\nu\in\Omega^\circ$, a basis $b_{\nu,1},\ldots,b_{\nu,d_\nu}$ of $V_\lambda[\nu]$, where $d_\nu:=\dim V_\lambda[\nu]$; thus, $\sum_{\nu\in\Omega^\circ}d_\nu=a(\lambda)$, since every weight of $V_\lambda$ is a weight of $\sX_1^\C$.
Since $P_\lambda$ commutes with $\operatorname{ad}_{R_X}$ for $X\in\frakt$, it maps $\sX_1^\C[\nu]$ onto $I_\lambda[\nu]=V_\lambda[\nu]\otimes C_\lambda$.
Hence, for every $w$ in the complexification of $\sX_1^\circ$ and $\nu\in\Omega^\circ$, there are unique vectors $c_{\nu,1}(w),\ldots,c_{\nu,d_\nu}(w)\in C_\lambda$ such that
\begin{equation}\label{eq:components}
P_\lambda w_\nu=\sum_{i=1}^{d_\nu}b_{\nu,i}\otimes c_{\nu,i}(w),
\end{equation}
and the map $w\mapsto(c_{\nu,i}(w))_{\nu,i}$ from the complexification of $\sX_1^\circ$ to $C_\lambda^{a(\lambda)}$ is linear and surjective, since the components $w_\nu$ range independently over the spaces $\sX_1^\C[\nu]$, for $\nu\in\Omega^\circ$.
Let $C'_\lambda(w)\subseteq C_\lambda$ be the span of the vectors $c_{\nu,i}(w)$, for $\nu\in\Omega^\circ$ and $i=1,\ldots,d_\nu$.
We claim that
\begin{equation}\label{eq:M-w}
\mathcal M_w=\bigoplus_\lambda V_\lambda\otimes C'_\lambda(w),
\end{equation}
where $\lambda$ ranges over the highest weights with $c_\lambda(\sX_1^\C)\neq0$; in particular, $\mathcal M_w=\sX_1^\C$ if and only if $C'_\lambda(w)=C_\lambda$ for all these $\lambda$.

To prove~\eqref{eq:M-w}, let $\nu\in\Omega^\circ$.
For $H\in\frakz$ and $X\in\fraks^\C$, we have that $[R_H,R_X]=0$, so $\operatorname{ad}_{R_H}$ commutes with the action of $\mathcal U(\fraks^\C)$ and acts on $\mathcal U(\fraks^\C)w_\nu$ by the scalar $\nu(H)$, as on $w_\nu$.
Thus, $\mathcal U(\fraks^\C)w_\nu$ is invariant under $\operatorname{ad}_{R_X}$ for all $X\in\frakk^\C=\frakz^\C \oplus \fraks^\C$ and hence, it is a $K$-submodule of $\sX_1^\C$ (note that $K$ is connected). Consequently, $\mathcal M_w$ is a $K$-submodule of $\sX_1^\C$.
Every $K$-submodule $N$ of $\sX_1^\C$ is completely reducible, and each of its irreducible submodules of type $V_\lambda$ is contained in $I_\lambda$; thus, $N=\bigoplus_\lambda(N\cap I_\lambda)$ and $P_\lambda(N)=N\cap I_\lambda$.
Applying this to $N=\mathcal M_w$, and using that $P_\lambda$ commutes with the action of $\mathcal U(\fraks^\C)$, we obtain that
$$
\mathcal M_w=\bigoplus_\lambda P_\lambda(\mathcal M_w)
\quad\mbox{and}\quad
P_\lambda(\mathcal M_w)=\sum_{\nu\in\Omega^\circ}\mathcal U(\fraks^\C)\,P_\lambda w_\nu.
$$
By~\eqref{eq:components}, each $P_\lambda w_\nu$ belongs to $V_\lambda\otimes C'_\lambda(w)$, which is invariant under $\fraks^\C$; thus, $P_\lambda(\mathcal M_w)\subseteq V_\lambda\otimes C'_\lambda(w)$.
Conversely, let $\nu\in\Omega^\circ$, $i_0\in\{1,\ldots,d_\nu\}$, and $u\in V_\lambda$.
Since $\frakz$ acts on the irreducible $K$-module $V_\lambda$ by scalars, $V_\lambda$ is irreducible as an $\fraks^\C$-module, and by Burnside's theorem (see, e.g.,~\cite[Chapter~XVII, Section~3]{GTM211}), every endomorphism of $V_\lambda$ is the action of an element of $\mathcal U(\fraks^\C)$.
We choose $a\in\mathcal U(\fraks^\C)$ that acts on $V_\lambda$ by a linear map $A$ with $Ab_{\nu,i_0}=u$ and $Ab_{\nu,i}=0$ for $i\neq i_0$.
Since $\mathcal U(\fraks^\C)$ acts on $V_\lambda\otimes C_\lambda$ through the first factor, \eqref{eq:components} gives $a\cdot P_\lambda w_\nu=\sum_iAb_{\nu,i}\otimes c_{\nu,i}(w)=u\otimes c_{\nu,i_0}(w)$.
Thus, $V_\lambda\otimes c_{\nu,i_0}(w)\subseteq\mathcal U(\fraks^\C)P_\lambda w_\nu\subseteq P_\lambda(\mathcal M_w)$ for all $\nu$ and $i_0$, so $V_\lambda\otimes C'_\lambda(w)\subseteq P_\lambda(\mathcal M_w)$, which proves~\eqref{eq:M-w}.

\xc{Step~3: the Zariski open set.}
Let $\mathcal G$ be the set of $w\in\sX_1^\circ$ such that $C'_\lambda(w)=C_\lambda$ for every highest weight $\lambda$ with $c_\lambda(\sX_1^\C)\neq0$.
By Steps~1 and~2, $\frakw^\C\supseteq\sX_0^{\fraks,\C}+\sX_1^\C=(\sX_0^{\fraks}+\sX_1)^\C$ for $w\in\mathcal G$, and taking the $\sigma$-fixed points gives $\frakw\supseteq\sX_0^{\fraks}+\sX_1$.
It remains to show that $\mathcal G$ is a nonempty Zariski open subset of $\sX_1^\circ$.
Fix $\lambda$ with $c_\lambda(\sX_1^\C)\neq0$ and a basis of $C_\lambda$, and let $A_\lambda(w)$ be the $c_\lambda(\sX_1^\C)\times a(\lambda)$ matrix whose columns are the coordinates of the vectors $c_{\nu,i}(w)$, for $\nu\in\Omega^\circ$ and $i=1,\ldots,d_\nu$.
By condition~(W), $a(\lambda)\geq c_\lambda(\sX_1^\C)$, and $C'_\lambda(w)=C_\lambda$ if and only if some minor of $A_\lambda(w)$ of size $c_\lambda(\sX_1^\C)$ is nonzero.
The entries of $A_\lambda(w)$ are linear functions of $w$; thus, these minors are polynomial functions $p_{\lambda,1},\ldots,p_{\lambda,s_\lambda}$ on the complexification of $\sX_1^\circ$, and $\mathcal G$ is the set of $w\in\sX_1^\circ$ such that, for each $\lambda$, at least one of $p_{\lambda,1}(w),\ldots,p_{\lambda,s_\lambda}(w)$ is nonzero.
The complement of $\mathcal G$ in $\sX_1^\circ$ is thus the common zero set of the real and imaginary parts of the products $\prod_\lambda p_{\lambda,i_\lambda}$, for all choices of the indices $i_\lambda$, which are real polynomial functions on $\sX_1^\circ$; hence, $\mathcal G$ is Zariski open.
Finally, since $w\mapsto(c_{\nu,i}(w))_{\nu,i}$ is surjective onto $C_\lambda^{a(\lambda)}$ and $a(\lambda)\geq c_\lambda(\sX_1^\C)$, there is $w$ in the complexification of $\sX_1^\circ$ for which the vectors $c_{\nu,i}(w)$ span $C_\lambda$, so some $p_{\lambda,i_\lambda}$ is not the zero polynomial.
The product $p:=\prod_\lambda p_{\lambda,i_\lambda}$ of these nonzero polynomials is a nonzero polynomial function on the complexification of $\sX_1^\circ$, and it does not vanish identically on $\sX_1^\circ$, since a polynomial function on the complexification of a real vector space that vanishes on the real vector space is zero.
Every $w\in\sX_1^\circ$ with $p(w)\neq0$ belongs to $\mathcal G$, so $\mathcal G$ is nonempty.
\end{proof}

Combining Lemma~\ref{lem:separation} with Theorem~\ref{thm:reductive}, we obtain the following more precise form of Theorem~\ref{thm:main} for the case where $S$ acts nontrivially on $M$:

\begin{theorem}\label{thm:two-generators}
Suppose that $K$ is connected, that $S$ acts nontrivially on $M$, and that $(\mathbb V,v)$ satisfies condition~(W).
Let $T\in\frakt$ and $E\in\fraks$ be given as in~\eqref{eq:separation} and~\eqref{eq:E}, and $\mathcal G\subseteq\sX_1^\circ$ be given by Lemma~\ref{lem:separation}.
Then, $\xalg(M)=\operatorname{Lie}\{R_T,R_E+w\}$ for every $w\in\mathcal G$.
\end{theorem}

\begin{proof}
By Lemma~\ref{lem:separation}, $\operatorname{Lie}\{R_T,R_E+w\}$ contains $\sX_0^{\fraks}+\sX_1$, and hence, it contains $\xalg(M)$ by Theorem~\ref{thm:reductive}. Conversely, $R_T$, $R_E$, and $w$ belong to $\xalg(M)$.
\end{proof}

By Lemmas~\ref{lem:spread} and~\ref{lem:W}, equivariant embeddings that satisfy condition~(W) exist, namely, those given by Lemma~\ref{lem:spread} with $R\geq R_0$.
For spheres, condition~(W) holds for the standard embedding, as shown in Remark~\ref{rem:W-sphere}, and explicit choices of $w$ are given in Section~\ref{sec:spheres}; we do not know explicit choices of $w$ for other compact homogeneous spaces in general.
The following example shows that condition~(W) does not hold for every equivariant embedding:

\begin{example}\label{rem:W-fails}
Let $K=\mathrm{SO}_3(\R)$ act on $M=\mathrm{SO}_3(\R)$ by left multiplication, and consider the equivariant embedding $(\R^{3\times3},I)$,
  where $K$ acts on $\R^{3\times3}$ by left multiplication.
Since $\mathrm{SO}_3(\R)$ spans $\R^{3\times3}$ and the vector fields $R_X$, for $X$ in a basis of $\frakk$, are linearly independent at every point, the map $\ell\otimes X\mapsto(\ell|_M)R_X$ from $(\R^{3\times3})^*\otimes\frakk$ to $\sX_1$ is an isomorphism.
As a $K$-module, $\R^{3\times3}$ is the sum of three copies of $\R^3$, whose complexification is the adjoint module $V_\alpha$, where $\alpha$ is the positive root.
Since $V_\alpha\otimes V_\alpha\approx V_0\oplus V_\alpha\oplus V_{2\alpha}$, we have that $c_\alpha(\sX_1^\C)=3$, whereas $a(\alpha)=\dim V_\alpha[0]=1$.
 Thus, for every $w\in\sX_1^\circ$, the subspace $C'_\alpha(w)$ in the proof of Lemma~\ref{lem:separation} has dimension at most $a(\alpha)=1<3=c_\alpha(\sX_1^\C)$, so $\mathcal M_w\neq\sX_1^\C$ by~\eqref{eq:M-w}.
Nevertheless, $\operatorname{Lie}(\sX_0+\sX_1)=\xalg(M)$ for this embedding by Theorem~\ref{thm:every-embedding}; condition~(W) is only a sufficient condition for the argument of Lemma~\ref{lem:separation},
  and we do not know whether $\xalg(M)=\operatorname{Lie}\{R_T,R_E+w\}$ for some $w\in\sX_1$ in this example.
\end{example}

\subsection{Tori}\label{ssec:torus-two}

We now consider the case where $S$ acts trivially on $M$.
With the notation of Subsection~\ref{ssec:tori}, $M=\T^m$ with angle coordinates $\theta_1,\ldots,\theta_m$, and $\xalg(M)=\xtrig(\T^m)$.
We use the indices $j$ modulo $m$, so that $\theta_{m+1}:=\theta_1$.
In the following result, the terms $\cos3\theta_{j+1}$ in $W_2$ couple the coordinate $\theta_j$ to $\theta_{j+1}$, and the frequency $3$ keeps the eigenvalues of $\operatorname{ad}_{W_1}$ on the summands of $W_2$ pairwise distinct:

\begin{proposition}\label{prop:torus-two}
Let $m\in\N$ and $W_1:=\sum_{j=1}^m4^{j-1}\partial_j$.  Further, let $W_2:=(\cos\theta_1+\cos2\theta_1)\,\partial_1$ for $m=1$, and
$$
W_2:= \sum_{j=1}^m \bigl(\cos\theta_j+\cos2\theta_j+\cos3\theta_{j+1}\bigr) \, \partial_j    \quad\mbox{for } m\geq2.
$$
Then, $\xtrig(\T^m)=\operatorname{Lie}\{W_1,W_2\}$.
\end{proposition}

\begin{proof}
Let $\frakw:=\operatorname{Lie}\{W_1,W_2\}$.
By~\eqref{eq:torus-bracket}, $[W_1,e_\chi \partial_j]=\mathrm{i}\,(\chi\cdot t)\,e_\chi \partial_j$, where $t:=(1,4,\ldots,4^{m-1})$.
Writing $\cos\vartheta=\frac12(e^{\mathrm{i}\vartheta}+e^{-\mathrm{i}\vartheta})$, we decompose $W_2$ into the eigenvectors
$$
\tfrac12 e^{\pm\mathrm{i}\theta_j}\partial_j, \quad \tfrac12 e^{\pm2\mathrm{i}\theta_j}\partial_j, \quad\mbox{and, for  } m\geq2,\quad\tfrac12 e^{\pm3\mathrm{i}\theta_{j+1}}\partial_j,
$$
for $j=1,\ldots,m$, of $\operatorname{ad}_{W_1}$, with eigenvalues $\pm\mathrm{i}\,4^{j-1}$, $\pm2\mathrm{i}\,4^{j-1}$, and $\pm3\mathrm{i}\,4^{j}$, respectively, where $4^m$ is replaced by $1$ for $j=m$.
These eigenvalues are pairwise distinct, since each of them is of the form $c\,\mathrm{i}\,4^{k-1}$ with $c\in\{\pm1,\pm2,\pm3\}$ and $k\in\{1,\ldots,m\}$, each pair $(c,k)$ occurs at most once, and distinct pairs give distinct numbers $c\,4^{k-1}$, because $|c|\leq3<4$.
By Lemma~\ref{lem:interpolation}, all these eigenvectors belong to $\frakw^\C$.

For $c,c'\in\Z$, we have that $[e^{\mathrm{i}c\theta_j}\partial_j,e^{\mathrm{i}c'\theta_j}\partial_j]=\mathrm{i}(c'-c)e^{\mathrm{i}(c+c')\theta_j}\partial_j$.
In particular, $[e^{\mathrm{i}\theta_j}\partial_j,e^{-\mathrm{i}\theta_j}\partial_j]=-2\mathrm{i}\,\partial_j$, so $\partial_j\in\frakw^\C$.
By induction on $|c|$, we obtain that $e^{\mathrm{i}c\theta_j}\partial_j\in\frakw^\C$ for all $c\in\Z$: if $c\geq2$, then $[e^{\mathrm{i}\theta_j}\partial_j,e^{\mathrm{i}c\theta_j}\partial_j]=\mathrm{i}(c-1)e^{\mathrm{i}(c+1)\theta_j}\partial_j$,
and similarly for $c\leq-2$ with $e^{-\mathrm{i}\theta_j}\partial_j$ in place of $e^{\mathrm{i}\theta_j}\partial_j$.
For variables $\theta_{j_1},\ldots,\theta_{j_s}$, let $\mathcal T(\theta_{j_1},\ldots,\theta_{j_s})$ be the space of complex trigonometric polynomials in these variables.
 We have shown that $\mathcal T(\theta_j)\partial_j\subseteq\frakw^\C$ for $j=1,\ldots,m$; for $m=1$, this means that $\frakw^\C\supseteq\xtrig(\T^1)^\C$.

Now, let $m\geq2$, and fix $j$.
For $f\in\mathcal T(\theta_j)$, we have that $[f\partial_j,e^{3\mathrm{i}\theta_{j+1}}\partial_j]=-e^{3\mathrm{i}\theta_{j+1}}f'\partial_j$, where $f'$ is the derivative of $f$, since $\partial_je^{3\mathrm{i}\theta_{j+1}}=0$. Since the derivatives $f'$ span the trigonometric polynomials in $\theta_j$ with zero mean,
  and since $e^{3\mathrm{i}\theta_{j+1}}\partial_j\in\frakw^\C$, we obtain that $e^{3\mathrm{i}\theta_{j+1}}\mathcal T(\theta_j)\partial_j\subseteq\frakw^\C$.
Next, for $g\in\mathcal T(\theta_{j+1})$ and $f\in\mathcal T(\theta_j)$,
$$
  [g\partial_{j+1},e^{3\mathrm{i}\theta_{j+1}}f\partial_j]=3\mathrm{i}\,g\,e^{3\mathrm{i}\theta_{j+1}}f\,\partial_j,
$$
since $\partial_jg=0$. As $g$ ranges over $\mathcal T(\theta_{j+1})$, so does $ge^{3\mathrm{i}\theta_{j+1}}$, and hence, $\mathcal T(\theta_j,\theta_{j+1})\partial_j\subseteq\frakw^\C$.
 For $m=2$, this means that $\frakw^\C\supseteq\xtrig(\T^2)^\C$.

For $m\geq3$, we show by induction on $s=2,\ldots,m$ that $\mathcal T(\theta_j,\theta_{j+1},\ldots,\theta_{j+s-1})\partial_j\subseteq \frakw^\C$ for all $j$.
The case $s=2$ has just been shown. Let $2\leq s\leq m-1$, $F\in\mathcal T(\theta_{j+1},\ldots,\theta_{j+s})$, and $f\in\mathcal T(\theta_j)$.
By the induction hypothesis for $j+1$ in place of $j$, $F\partial_{j+1}\in \frakw^\C$, and by the case $s=2$, $e^{\mathrm{i}\theta_{j+1}} f \partial_j\in \frakw^\C$.
Since $s\leq m-1$, the function $F$ does not depend on $\theta_j$, and hence,
$$
[F\partial_{j+1},e^{\mathrm{i}\theta_{j+1}}f\partial_j] = \mathrm{i}\,Fe^{\mathrm{i}\theta_{j+1}}f\,\partial_j.
$$
As $F$ ranges over $\mathcal T(\theta_{j+1},\ldots,\theta_{j+s})$,
  so does $Fe^{\mathrm{i}\theta_{j+1}}$; thus, $\mathcal T(\theta_j,\ldots,\theta_{j+s})\partial_j\subseteq \frakw^\C$.
 For $s=m$, we obtain that $\frakw^\C$ contains $\mathcal T(\theta_1,\ldots,\theta_m)\partial_j$ for all $j$, i.e., $\frakw^\C\supseteq \xtrig(\T^m)^\C$ also for $m\geq3$.
In all cases, taking the $\sigma$-fixed points gives $\xtrig(\T^m)\subseteq \frakw$, and the converse inclusion holds since $W_1$ and $W_2$ are trigonometric.
\end{proof}

Previously, $\xtrig(\T^m)$, for $m\geq2$, was known to be generated by the $m+2$ explicit vector fields of~\cite[Theorem~1]{Luk79} and by the $4m$ explicit vector fields of~\cite[Lemma~6.5]{AS22}; for $m=2,3$, the divergence-free trigonometric vector fields are generated by $\partial_1,\ldots,\partial_m$ together with one divergence-free trigonometric vector field whose frequencies generate $\Z^m$ and satisfy a nondegeneracy condition~\cite[Theorems~6 and~8]{AK26}.

We are now in a position to prove Theorem~\ref{thm:main}:

\begin{proof}[Proof of Theorem~\ref{thm:main}]
By Lemma~\ref{lem:subgroup}, we can assume that $K$ is connected.
If $S$ acts nontrivially on $M$, then Lemmas~\ref{lem:spread} and~\ref{lem:W} provide an equivariant embedding that satisfies condition~(W), and Theorem~\ref{thm:two-generators} gives $\xalg(M)=\operatorname{Lie}\{R_T,R_E+w\}$ for every $w\in\mathcal G$.
If $S$ acts trivially on $M$, then by Subsection~\ref{ssec:tori}, $M$ is identified with $\T^m$,
    $m=\dim M$, so that $K$ acts by translations and $\xalg(M)=\xtrig(\T^m)$, and Proposition~\ref{prop:torus-two} applies.
\end{proof}

\section{Spheres}\label{sec:spheres}

This section has two purposes: to show that Theorems~\ref{thm:main}, \ref{thm:every-embedding}, and~\ref{thm:two-generators} apply to spheres with the standard embedding, and to obtain explicit generators from the results of~\cite{Luk77}. 
In this section, we let $n\geq2$ and consider the sphere $S^{n-1}=\{x\in\R^n\mid\sum_{j=1}^nx_j^2=1\}$ with the action of $K=\mathrm{SO}_n(\R)$ by matrix multiplication.
Let $e_1,\ldots,e_n$ be the standard basis of $\R^n$; then, $(\R^n,e_n)$ is an equivariant embedding of $S^{n-1}$.
We extend the Euclidean inner product to the bilinear form $\langle a,b\rangle:=a^{\mathsf T}b$ on $\C^n$, and for $a,b\in\C^n$, we let
$$
\ell_a(x):=\langle a,x\rangle \quad\mbox{and}\quad
a\wedge b:=ab^{\mathsf T}-ba^{\mathsf T}\in\mathfrak{so}_n(\C).
$$
For $A\in\mathfrak{so}_n(\C)$, we have that $R_A(x)=-Ax$, and in particular, $R_{a\wedge b}(x)=\ell_a(x)b-\ell_b(x)a$.
For $1\leq i<j\leq n$, we let $R_{ij}:=R_{e_i\wedge e_j}=x_i\partial/\partial x_j-x_j\partial/\partial x_i$ be the infinitesimal rotation in the $x_ix_j$-plane.

By Lemma~\ref{lem:kfinite}, applied to $(\R^n,e_n)$, the algebra $\Ocal(S^{n-1})$ comprises the restrictions to $S^{n-1}$ of polynomials in $x_1,\ldots,x_n$, and by Lemma~\ref{lem:tangent}, $\xalg(S^{n-1})$ is spanned by the vector fields $fR_{ij}$ with $f\in\Ocal(S^{n-1})$.
Thus, $\xalg(S^{n-1})$ comprises the restrictions to $S^{n-1}$ of the polynomial vector fields on $\R^n$ that are tangent to $S^{n-1}$.
To wit, let $\sum_{k=1}^nd_k\,\partial/\partial x_k$ be a polynomial vector field with $\sum_{k=1}^nx_kd_k=0$ on $S^{n-1}$. Then, it coincides on $S^{n-1}$ with $\sum_{i<j}(x_id_j-x_jd_i)R_{ij}$, since the $k$th components of both vector fields are equal to $d_k$ on $S^{n-1}$.
Thus, $\xalg(S^{n-1})$ is the Lie algebra of algebraic vector fields on $S^2$ considered in~\cite{CofDS2} for $n=3$,
and it is the Lie algebra of spherical vector fields considered in~\cite{Luk77}, which are defined as the smooth vector fields whose translates under $\mathrm{SO}_n(\R)$ span a finite-dimensional space.

By Theorem~\ref{thm:main}, $\xalg(S^{n-1})$ is generated by two vector fields for every $n\geq2$.
By Theorem~\ref{thm:every-embedding}, applied to $(\R^n,e_n)$, we have that $\xalg(S^{n-1})=\operatorname{Lie}(\sX_0+\sX_1)$ for $n\geq3$; for $n=2$, the circle $S^1\subseteq\R^2$ is the exceptional case of Theorem~\ref{thm:every-embedding}, since $\mathcal A_1=\operatorname{span}\{\cos\theta,\sin\theta\}$ for the angle coordinate $\theta$, and two generators of $\xalg(S^1)$ are given in Proposition~\ref{prop:torus-two}.
In the remainder of the section, we let $n\geq3$.

\subsection{Explicit generators}\label{ssec:lukatskii}

We now show that two explicit generators of the divergence-free part of $\xalg(S^{n-1})$ and of $\xalg(S^{n-1})$ itself are obtained, for every $n\geq3$, by combining two results of~\cite{Luk77} with Lemma~\ref{lem:interpolation}; the case $n=3$ of the second result is proved in Subsection~\ref{ssec:S2}.
Let $r:=\lfloor n/2\rfloor$, and for $j=1,\ldots,r$, let
$$
    u_j:=e_{2j-1}+\mathrm{i}e_{2j} \quad\mbox{and}\quad
    J_j:=e_{2j-1}\wedge e_{2j}.
$$
Then, $\frakh:=\operatorname{span}_\R\{J_1,\ldots,J_r\}$ is a maximal abelian subalgebra of $\mathfrak{so}_n(\R)$, and since $\mathfrak{so}_n(\R)$ is semi-simple, we can take $\frakt=\frakh$ in Subsection~\ref{ssec:weights}, whose notation for weights and roots we use; in particular, $r=\dim\frakh=\lfloor n/2\rfloor$ agrees with the notation there.
Let $\varepsilon_j\in(\frakh^\C)^*$ be given by $\varepsilon_j(J_k)=\mathrm{i}\delta_{jk}$.
Then, $J_ku_j=\mathrm{i}\delta_{jk}u_j$ and $J_k\sigma(u_j)=-\mathrm{i}\delta_{jk}\sigma(u_j)$, and for odd $n$, $J_ke_n=0$; thus, $u_j$, $\sigma(u_j)$, and $e_n$ (for odd $n$) are weight vectors in $\C^n$ of weights $\varepsilon_j$, $-\varepsilon_j$, and $0$, respectively.
Since $[X,a\wedge b]=(Xa)\wedge b+a\wedge(Xb)$ for $X\in\mathfrak{so}_n(\C)$, the roots of $\mathfrak{so}_n(\C)$ with respect to $\frakh^\C$ are $\pm\varepsilon_i\pm\varepsilon_j$, for $1\leq i<j\leq r$, and, for odd $n$, $\pm\varepsilon_j$, for $1\leq j\leq r$.
For instance, the root spaces of $\varepsilon_i-\varepsilon_j$, $\varepsilon_i+\varepsilon_j$, and $\varepsilon_i$ are spanned by $u_i\wedge\sigma(u_j)$, $u_i\wedge u_j$, and $u_i\wedge e_n$, respectively.
We choose the simple roots $\varepsilon_j-\varepsilon_{j+1}$, for $1\leq j<r$, together with $\varepsilon_r$ for odd $n$ and with $\varepsilon_{r-1}+\varepsilon_r$ for even $n$, and let
\begin{equation}\label{eq:sphere-E}
E_+ := \sum_{j=1}^{r-1}u_j\wedge\sigma(u_{j+1})+
\begin{cases}
u_r \wedge e_n, & n\mbox{ odd},\\
u_{r-1} \wedge u_r, & n\mbox{ even},
\end{cases}
\quad \mbox{and} \quad
E:= E_+ + \sigma(E_+)\in\mathfrak{so}_n(\R),
\end{equation}
so that $E_+$ is the sum of nonzero simple root vectors, as in Lemma~\ref{lem:separation}.
For $X\in\mathfrak{so}_n(\C)$, $a\in\C^n$, and $B\in\mathfrak{so}_n(\C)$, we have that $R_X\ell_a=\ell_{Xa}$, and hence, $[R_X,\ell_aR_B]=\ell_{Xa}R_B+\ell_aR_{[X,B]}$.

For $a\in\C^n$, let $G_a(x):= a-\ell_a(x)x$; for $a\in\R^n$, the vector field $G_a$ is the gradient of $\ell_a$ on $S^{n-1}$, and we write $G_{x_j}:=G_{e_j}$.
A direct computation shows that
\begin{equation}\label{eq:conformal}
[R_X,G_a]=G_{Xa} \quad\mbox{and}\quad
[G_a,G_b]=R_{a\wedge b}, \quad\mbox{for } X\in\mathfrak{so}_n(\C)\mbox{ and }a,b\in\C^n.
\end{equation}
Thus, the space $\frakc:=\sX_0\oplus\{G_a \mid  a\in\R^n\}$ of conformal vector fields on $S^{n-1}$ is a Lie subalgebra of $\xalg(S^{n-1})$; the sum is direct since $R_A(-x)=-R_A(x)$ and $G_a(-x)=G_a(x)$, and since $a\mapsto G_a$ is injective.
 Let $\omega$ be the volume form of $S^{n-1}$ induced by the Euclidean metric, $\operatorname{div}V$ be the divergence of a smooth vector field $V$ with respect to $\omega$, and $\operatorname{Diff}_0(S^{n-1},\omega)$ be the group of $\omega$-preserving diffeomorphisms of $S^{n-1}$ that are smoothly isotopic to $\operatorname{id}_{S^{n-1}}$ through $\omega$-preserving diffeomorphisms.
Since $\operatorname{div}[U,V]=U(\operatorname{div}V)-V(\operatorname{div}U)$, the divergence-free algebraic vector fields form a Lie subalgebra $\xalgz(S^{n-1})$ of $\xalg(S^{n-1})$.

The following theorem is the Lie-algebraic part of~\cite[Theorem~1]{Luk77}; it was announced in~\cite[Theorem~3]{Luk74}, with $x_nR_{12}$ in place of $x_3R_{12}$.
In~\cite{Luk77}, it is stated with the rotation fields $R_{i,i+1}$, for $i=1,\ldots,n-1$, in place of $\sX_0$; these fields generate $\sX_0$.
It is also stated there for a Lie algebra that is identified with $\xalgz(S^{n-1})$ in~\cite[Lemma~1]{Luk77}.
The remaining assertions of~\cite[Theorem~1]{Luk77} concern the topological Lie algebra of divergence-free smooth vector fields on $S^{n-1}$ and the group of volume-preserving diffeomorphisms of $S^{n-1}$, as well as real projective spaces; the assertion on the group is recalled in the discussion of Table~\ref{tab:counts}.
The proof in~\cite{Luk77} applies to every $n\geq3$, including $n=3$, since it is carried out inside $\xalgz(S^{n-1})$, whose irreducible $\mathrm{SO}_n(\R)$-summands have pairwise distinct highest weights for every $n\geq3$ (see the decomposition in~\cite[p.~120]{Luk77}); the isomorphic pairs of summands discussed before Theorem~\ref{thm:luk2} below involve summands of $\xalg(S^2)$ that lie outside $\xalgz(S^2)$.

\begin{theorem}[cf.~{\cite[Theorem~1]{Luk77}}]\label{thm:luk1}
For $n\geq3$, it holds that $\xalgz(S^{n-1})=\operatorname{Lie}(\sX_0\cup\{x_3R_{12}\})$.
\end{theorem}

The following theorem is stated in~\cite[Lemma~2]{Luk77} for spheres of all dimensions.
The proof there uses the decomposition of $\xalg(S^{n-1})$ into irreducible $\mathrm{SO}_n(\R)$-modules and the fact that a Lie subalgebra of $\xalg(S^{n-1})$ that contains $\sX_0$ and a vector field with a nonzero component in an irreducible summand contains the entire summand.
This fact holds for $n\geq4$, since the irreducible summands are then pairwise nonisomorphic (see the decomposition in~\cite[p.~120]{Luk77}).
For $n=3$, however, the irreducible summands occur in isomorphic pairs, and a nonzero invariant subspace of the sum of such a pair need not contain either summand, as the graph of an isomorphism between them shows.
Thus, for $n=3$, this fact does not follow from irreducibility, and the proof requires an additional argument; this case is not treated separately in that proof.
In Subsection~\ref{ssec:S2}, we prove the case $n=3$ by an argument that uses the Casimir operator of $\mathfrak{so}_3(\R)$, the conformality of the flows of the vector fields $G_a$, and the fact that the irreducible $\mathrm{SO}_3(\R)$-modules are odd-dimensional.

\begin{theorem}\label{thm:luk2}
For $n\geq3$, the Lie algebra $\frakc$ is a maximal Lie subalgebra of $\xalg(S^{n-1})$, i.e.,
every Lie subalgebra of $\xalg(S^{n-1})$ that strictly contains $\frakc$ is $\xalg(S^{n-1})$.
\end{theorem}

For $n\geq4$, Theorem~\ref{thm:luk2} is~\cite[Lemma~2]{Luk77}; the proof for $n=3$ is given in Subsection~\ref{ssec:S2}.

For $n=3$, we also use the notation of~\cite{CofDS2}: for a polynomial $f$ in $x_1,x_2,x_3$, let
$$
  X_f:=x\times\nabla f
  \quad\mbox{and}\quad
  G_f:=\nabla f-(x\cdot\nabla f)x,
$$
restricted to $S^2$, where $\times$ is the cross product on $\R^3$.
Then, $G_{\ell_a}=G_a$ for $a\in\R^3$, and since $x\times e_1=-R_{23}(x)$, $x\times e_2=R_{13}(x)$, and $x\times e_3=-R_{12}(x)$, we have that
$$
X_f=-(\partial_1f)R_{23}+(\partial_2f)R_{13}-(\partial_3f)R_{12},
$$
where $\partial_j:=\partial/\partial x_j$; in particular, $X_{x_1}=-R_{23}$, $X_{x_2}=R_{13}$, $X_{x_3}=-R_{12}$, and $X_f\in\sX_1$ if $f$ is a homogeneous polynomial of degree two.
By~\cite[Section~2.2]{CofDS2}, rotating $f$ rotates $X_f$, so $[R_X,X_f]=X_{R_Xf}$ for $X\in\mathfrak{so}_3(\C)$ by~\eqref{eq:R-derivatives}, where $X_f$ is defined for complex polynomials $f$ by complex linearity.

Let $\Sigma\subseteq(\frakh^\C)^*$ be the set consisting of the simple roots, their negatives, $\pm\varepsilon_1$, and $0$ for odd $n$ or $\pm\varepsilon_r$ for even $n$, and let $T\in\frakh$ be such that
\begin{equation}\label{eq:sphere-separation}
\lambda(T)\neq\lambda'(T)\quad\mbox{whenever }\lambda,\lambda'\in\Sigma\mbox{ are distinct};
\end{equation}
for example, $T:=\sum_{j=1}^r4^{r-j}J_j$ satisfies~\eqref{eq:sphere-separation}, since $\varepsilon_j(T)=\mathrm{i}\,4^{r-j}$, so that the values at $T$ of the elements of $\Sigma$ are $0$ or $\pm\mathrm{i}$ (for odd $n$) or $\pm5\mathrm{i}$ (for even $n$), $\pm3\mathrm{i}\,4^{r-j-1}$ for $1\leq j<r$, and $\pm\mathrm{i}\,4^{r-1}$, which are pairwise distinct.
 Combining Theorems~\ref{thm:luk1} and~\ref{thm:luk2} with Lemma~\ref{lem:interpolation}, we obtain the following explicit generators:

\begin{proposition}\label{prop:sphere-pairs}
Let $n\geq3$, and let $E$ and $T$ be given as in~\eqref{eq:sphere-E} and~\eqref{eq:sphere-separation}. Then, the following hold:
\begin{enumerate}
\item $\xalgz(S^{n-1})=\operatorname{Lie}\{R_T,\,R_E+x_nR_{12}\}$;
\item if $n\geq4$, then $\xalg(S^{n-1})=\operatorname{Lie}\{R_T,\,R_E+x_nR_{12}+G_{x_1}\}$;
	\item if $n=3$, then $\xalg(S^2)=\operatorname{Lie}\{R_T,\,R_E+X_{x_1^2-x_2^2}+G_{x_3}\}$.
\end{enumerate}
\end{proposition}

\begin{proof}
Let $W_2':=R_E+x_nR_{12}$ and $W_2:=W_2'+G_{x_1}$, and let $E_\alpha$ be the summand of $E_+$ of weight $\alpha$, for each simple root $\alpha$.
Then, $R_{E_\alpha}$ and $R_{\sigma(E_\alpha)}$ have weights $\alpha$ and $-\alpha$, respectively, and $G_{x_1}=\frac12(G_{u_1}+G_{\sigma(u_1)})$, where $G_{u_1}$ and $G_{\sigma(u_1)}$ have weights $\varepsilon_1$ and $-\varepsilon_1$ by~\eqref{eq:conformal}.
Since $[R_X,R_{12}]=0$ for $X\in\frakh^\C$, the vector field $\ell_aR_{12}$ has weight $\mu$ if $a$ is a weight vector of weight $\mu$.
Thus, for odd $n$, the vector field $x_nR_{12}=\ell_{e_n}R_{12}$ has weight $0$.
For even $n$, since $e_n=\frac{1}{2\mathrm{i}}(u_r-\sigma(u_r))$, it is the sum of the vector fields $\frac{1}{2\mathrm{i}}\ell_{u_r}R_{12}$ and $-\frac{1}{2\mathrm{i}}\ell_{\sigma(u_r)}R_{12}$, which have weights $\varepsilon_r$ and $-\varepsilon_r$.
It follows that $W_2'$ and, for $n\geq4$, $W_2$ are sums of nonzero vector fields whose weights are pairwise distinct elements of $\Sigma$: the weights $0$ and $\pm\varepsilon_r$ of the components of $x_nR_{12}$,
  for odd and even $n$, respectively, are not roots, and for $n\geq4$, the weights $\pm\varepsilon_1$ are not simple roots or their negatives, and they differ from $\pm\varepsilon_r$ since $r\geq2$.
These vector fields are eigenvectors of $\operatorname{ad}_{R_T}$, whose eigenvalues are pairwise distinct by~\eqref{eq:sphere-separation}.

Let $\frakw$ be either $\operatorname{Lie}\{R_T,W_2'\}$ or, for $n\geq4$, $\operatorname{Lie}\{R_T,W_2\}$.
By Lemma~\ref{lem:interpolation}, the summands of $W_2'$, resp., of $W_2$, listed above belong to $\frakw^\C$.
In particular, $\frakw^\C$ contains $R_{E_\alpha}$ and $R_{\sigma(E_\alpha)}$ for every simple root $\alpha$, and since these root vectors generate $\mathfrak{so}_n(\C)$, we have that $\sX_0\subseteq\frakw$ after taking the $\sigma$-fixed points.
Hence, $x_nR_{12}$, which is the sum of the summands listed above of weight $0$ for odd $n$, or of weights $\pm\varepsilon_r$ for even $n$, belongs to $\frakw^\C$, and hence, to $\frakw$, since it is real.
Since $\frakw$ contains $\sX_0$, it is invariant under the action of $\mathrm{SO}_n(\R)$.
To wit, let $V\in\frakw$ and $A\in\mathfrak{so}_n(\R)$. The translates of $V$ span a finite-dimensional space that is invariant under $\operatorname{ad}_{R_A}$, in which $t\mapsto e^{tA}\cdot V$ solves the linear differential equation $\dot Y=\operatorname{ad}_{R_A}Y$ by~\eqref{eq:R-derivatives}. Thus, $e^{tA}\cdot V=e^{t\operatorname{ad}_{R_A}}V\in\frakw$, and such elements $e^{tA}$ generate $\mathrm{SO}_n(\R)$.
For $n\geq4$, let $g\in\mathrm{SO}_n(\R)$ fix $e_1$ and $e_2$ and map $e_n$ to $e_3$. Then,
  $g\cdot x_n=x_3$ and $g\cdot R_{12}=R_{ge_1\wedge ge_2}=R_{12}$, so $x_3R_{12}=g\cdot(x_nR_{12})\in\frakw$; for $n=3$, this holds trivially.

Now, let $\frakw:=\operatorname{Lie}\{R_T,W_2'\}$.
By Theorem~\ref{thm:luk1}, $\xalgz(S^{n-1})\subseteq\frakw$.
Conversely, $R_T$ and $R_E$ are divergence-free, since their flows consist of rotations, and so is $x_nR_{12}$, since $\operatorname{div}(x_nR_{12})=R_{12}x_n=0$; thus, $\frakw\subseteq\xalgz(S^{n-1})$.

Next, let $n\geq4$ and $\frakw:=\operatorname{Lie}\{R_T,W_2\}$.
By~\eqref{eq:conformal}, the linear map $a\mapsto G_a$ from $\C^n$ to $\xalg^\C(S^{n-1})$ is an injective homomorphism of $\mathfrak{so}_n(\C)$-modules, and $\C^n$ is an irreducible $\mathfrak{so}_n(\C)$-module. Since $\frakw^\C$ contains $\sX_0^\C$ and $G_{u_1}$, it contains $\{G_a\mid a\in\C^n\}$, and hence, $\frakc\subseteq\frakw$ after taking the $\sigma$-fixed points.
The vector field $x_nR_{12}\in\frakw$ does not belong to $\frakc$.
To wit, suppose that $x_nR_{12}=R_A+G_a$ with $A\in\mathfrak{so}_n(\R)$ and $a\in\R^n$.  Replacing $x$ by $-x$ gives $x_nR_{12}=-R_A+G_a$, so $x_nR_{12}=G_a$. Evaluating at $e_n$ and at $e_1$ gives $a-a_ne_n=0$ and $a-a_1e_1=0$, so $a=0$, which contradicts $x_nR_{12}\neq0$.
By Theorem~\ref{thm:luk2}, $\frakw=\xalg(S^{n-1})$.

Finally, let $n=3$, $W_2'':=R_E+X_{x_1^2-x_2^2}+G_{x_3}$, and $\frakw:=\operatorname{Lie}\{R_T,W_2''\}$.
Then, $r=1$, the only simple root is $\varepsilon_1$, $\Sigma=\{0,\pm\varepsilon_1\}$, and~\eqref{eq:sphere-separation} means that $T$ is a nonzero element of $\frakh=\R J_1$; thus, the values of $0$,
  $\pm\varepsilon_1$, and $\pm2\varepsilon_1$ at $T$ are pairwise distinct.
Since $x_1^2-x_2^2=\frac12(\ell_{u_1}^2+\ell_{\sigma(u_1)}^2)$ and $R_X\ell_{u_1}^2=2\ell_{u_1}\ell_{Xu_1}$ for $X\in\mathfrak{so}_3(\C)$, the functions $\ell_{u_1}^2$ and $\ell_{\sigma(u_1)}^2$ have weights $2\varepsilon_1$ and $-2\varepsilon_1$,
    and so do the vector fields $X_{\ell_{u_1}^2}$ and $X_{\ell_{\sigma(u_1)}^2}$, whose sum is $2X_{x_1^2-x_2^2}$.
Thus, $W_2''$ is the sum of the nonzero vector fields $R_{E_\alpha}$, $R_{\sigma(E_\alpha)}$, $\frac12X_{\ell_{u_1}^2}$, $\frac12X_{\ell_{\sigma(u_1)}^2}$, and $G_{e_3}$, whose weights $\varepsilon_1$, $-\varepsilon_1$, $2\varepsilon_1$, $-2\varepsilon_1$, and $0$ take pairwise distinct values at $T$, and by Lemma~\ref{lem:interpolation}, these vector fields belong to $\frakw^\C$.
As above, it follows that $\sX_0\subseteq\frakw$, and since $\frakw^\C$ contains $\sX_0^\C$ and $G_{e_3}$, the argument given above for $n\geq4$ shows that $\frakc\subseteq\frakw$.
The vector field $X_{x_1^2-x_2^2}=\frac12(X_{\ell_{u_1}^2}+X_{\ell_{\sigma(u_1)}^2})$ belongs to $\frakw$, and it does not belong to $\frakc$.
To wit, $X_{x_1^2-x_2^2}=-2x_1R_{23}-2x_2R_{13}$ is nonzero, it is invariant under $x\mapsto-x$, and it vanishes at $e_3$ and at $e_1$; hence, the argument given above for $x_nR_{12}$ shows that it is not of the form $R_A+G_a$.
By Theorem~\ref{thm:luk2}, $\frakw=\xalg(S^2)$.
\end{proof}

\begin{corollary}\label{cor:sphere-flows}
Let $n\geq3$, let $E$ and $T$ be given as in~\eqref{eq:sphere-E} and~\eqref{eq:sphere-separation}, and let $W_2:=R_E+x_nR_{12}$.
Then, the subgroup generated by the flows of $R_T$ and $W_2$ is dense in $\operatorname{Diff}_0(S^{n-1},\omega)$ with respect to the $\mathrm{C}^\infty$-topology, and there is $a>0$ such that $\operatorname{Diff}_0(S^{n-1},\omega)$ is topologically generated by $e^{aR_T}$, $e^{W_2}$, and $e^{\sqrt2W_2}$.
\end{corollary}

\begin{proof}
Let $\mathfrak X(S^{n-1},\omega)$ be the space of divergence-free smooth vector fields on $S^{n-1}$.
Since $\omega$ is invariant under $\mathrm{SO}_n(\R)$, we have that $\operatorname{div}(g\cdot V)=g\cdot\operatorname{div}V$ for $g\in\mathrm{SO}_n(\R)$ and $V\in\mathfrak X(S^{n-1})$; thus, the vector fields $V_\varphi$ in the proof of Corollary~\ref{cor:B} are divergence-free if $V$ is, and that proof shows that $\xalgz(S^{n-1})$ is dense in $\mathfrak X(S^{n-1},\omega)$.
By item~1 of Proposition~\ref{prop:sphere-pairs}, the closure of $\operatorname{Lie}\{R_T,W_2\}$ therefore contains $\mathfrak X(S^{n-1},\omega)$.
If $(P_t)_{t\in[0,1]}$ is a smooth isotopy from $\operatorname{id}_{S^{n-1}}$ to $f\in\operatorname{Diff}_0(S^{n-1},\omega)$ through $\omega$-preserving diffeomorphisms, then $V_t:=\dot P_t\circ P_t^{-1}$ is divergence-free, and the proof of Corollary~\ref{cor:B} shows that $f$ belongs to the closure of the subgroup generated by the flows of $R_T$ and $W_2$.
This subgroup is contained in $\operatorname{Diff}_0(S^{n-1},\omega)$, since $R_T$ and $W_2$ are divergence-free: the flows of $R_T$ and $R_E$ consist of rotations, and $\operatorname{div}(x_nR_{12})=R_{12}x_n=0$.
The second assertion follows as in the proof of Corollary~\ref{cor:three}.
\end{proof}

Two other explicit pairs of vector fields whose flows generate dense subgroups of $\operatorname{Diff}_0(S^{n-1})$ and of $\operatorname{Diff}_0(S^{n-1},\omega)$, respectively, have been given in~\cite[Theorem~1]{DeTurck}.

\begin{remark}\label{rem:sphere-form}
The generators in items~2 and~3 of Proposition~\ref{prop:sphere-pairs} are of the form $R_T$ and $R_E+w$ of Theorem~\ref{thm:two-generators}, with the explicit vector fields $w=x_nR_{12}+G_{x_1}$ and $w=X_{x_1^2-x_2^2}+G_{x_3}$ in $\sX_1$; the proof of Proposition~\ref{prop:sphere-pairs}, however, uses Theorem~\ref{thm:luk2} instead of Theorem~\ref{thm:reductive}, and it does not use condition~(W).
For $n=3$, the two generators $X_{x_3}$ and $X_{x_1}+X_{x_1^2-x_2^2}+G_{x_3}$ of $\xalg(S^2)$ given in~\cite{CofDS2} are those of Proposition~\ref{prop:sphere-pairs}(3) with $T=-J_1$ and with the root vector $\frac{\mathrm{i}}{2}u_1\wedge e_3$ in place of $u_1\wedge e_3$ in~\eqref{eq:sphere-E}.
To wit, $X_{x_3}=-R_{12}=R_{-J_1}$, and $X_{x_1}=-R_{23}=R_E$ for $E=E_\alpha+\sigma(E_\alpha)$ with $E_\alpha:=\frac{\mathrm{i}}{2}u_1\wedge e_3$, since $E_\alpha+\sigma(E_\alpha)=\frac{\mathrm{i}}{2}(u_1-\sigma(u_1))\wedge e_3=-e_2\wedge e_3$.
Thus, the proof of Proposition~\ref{prop:sphere-pairs}(3) applies to these two vector fields and gives another proof of the main result of~\cite{CofDS2}: there, the Lie algebra that they generate is shown by explicit brackets to contain $X_f$ and $G_f$ for every harmonic polynomial $f$, whereas here, by Theorem~\ref{thm:luk2}, it suffices that this Lie algebra strictly contains~$\frakc$.
\end{remark}

\begin{remark}\label{rem:W-sphere}
Condition~(W) holds for the equivariant embedding $(\R^n,e_n)$ of $S^{n-1}$, $n\geq3$, so that Theorem~\ref{thm:two-generators} applies to the standard embedding: with the notation of Section~\ref{sec:two-generators}, $\xalg(S^{n-1})=\operatorname{Lie}\{R_T,R_E+w\}$ for every $T\in\frakh$ satisfying~\eqref{eq:separation} and every $w$ in a nonempty Zariski open subset of $\sX_1^\circ$.
To wit, $((\R^n)^*)^\C\approx\C^n=V_{\varepsilon_1}$, so $\Lambda(\R^n)=\{\varepsilon_1\}$, and~\eqref{eq:multiplicity-bound} gives $c_\lambda(\sX_1^\C)\leq\dim\mathfrak{so}_n(\C)[\lambda-\varepsilon_1]$ for every highest weight $\lambda$; thus,
    $c_{\varepsilon_1}(\sX_1^\C)\leq r$ and $c_\lambda(\sX_1^\C)\leq1$ for $\lambda\neq\varepsilon_1$.
Among the $n$ weights $\pm\varepsilon_j$ and, for odd $n$, $0$ of $\C^n$, counted with multiplicity, at most two, namely, $\pm\varepsilon_r$ for odd $n$, belong to $\Pi\cup(-\Pi)$; thus, $a(\varepsilon_1)\geq n-2\geq r$.
For every highest weight $\lambda$, we have that $a(\lambda)\geq1$: if $\lambda\notin\Pi\cup(-\Pi)$, then $\lambda$ itself is a weight of $V_\lambda$ outside $\Pi\cup(-\Pi)$; the negative of a simple root is not dominant; and the only dominant simple roots are $\varepsilon_1$ for $n=3$ and $\varepsilon_1\pm\varepsilon_2$ for $n=4$, in which cases $0$ is a weight of $V_\lambda$ outside $\Pi\cup(-\Pi)$.
Hence, $c_\lambda(\sX_1^\C)\leq a(\lambda)$ for every highest weight $\lambda$.
\end{remark}

 Table~\ref{tab:counts} lists the numbers of vector fields that generate $\xalgz(S^{n-1})$ and $\xalg(S^{n-1})$ as Lie algebras, and the numbers of flows and of diffeomorphisms that topologically generate $\operatorname{Diff}_0(S^{n-1},\omega)$ and $\operatorname{Diff}_0(S^{n-1})$, i.e., that generate a dense subgroup.
The first row is obtained from~\cite{Luk77}.
By Theorem~\ref{thm:luk1}, the $n$ vector fields $R_{12},\ldots,R_{n-1,n},x_3R_{12}$ generate $\xalgz(S^{n-1})$, since $R_{12},\ldots,R_{n-1,n}$ generate $\sX_0$.
Together with $G_{x_n}$, they generate $\xalg(S^{n-1})$ by Theorem~\ref{thm:luk2}, since the Lie algebra generated by them contains $\frakc$ by~\eqref{eq:conformal}, as well as $x_3R_{12}\notin\frakc$.
By~\cite[Theorem~1]{Luk77}, $\operatorname{Diff}_0(S^{n-1},\omega)$ is topologically generated by four diffeomorphisms, namely, the time-one maps of two vector fields in $\sX_0$ and the maps at times $1$ and $\sqrt2$ of the flow of $x_3R_{12}$; they lie on three flows.
By~\cite[Theorem~3]{Luk77}, $\operatorname{Diff}_0(S^{n-1})$ is topologically generated by two diffeomorphisms, namely, a rotation and the composition of a rotation with fixed-time maps of the flows of $G_{x_n}$ and $(x_1^2+x_2^2)R_{12}$; they lie in the subgroup generated by four flows.
For $n=3$, the entries of the first row for $\xalg(S^2)$ and $\operatorname{Diff}_0(S^2)$ rely on~\cite[Lemma~2]{Luk77} for $S^2$, i.e., on the case $n=3$ of Theorem~\ref{thm:luk2}, which is proved in Subsection~\ref{ssec:S2}, and the fixed time of the flow of $(x_1^2+x_2^2)R_{12}$ in~\cite[Theorem~3]{Luk77} is then restricted to sufficiently small nonzero values.

The second row is obtained from the pairs of vector fields in Proposition~\ref{prop:sphere-pairs}, which combine Theorems~\ref{thm:luk1} and~\ref{thm:luk2} with Lemma~\ref{lem:interpolation}, the interpolation argument in the proof of Kuranishi's theorem; their flows topologically generate the groups by Corollaries~\ref{cor:B} and~\ref{cor:sphere-flows}, and three diffeomorphisms are given by Corollaries~\ref{cor:three} and~\ref{cor:sphere-flows}.
The third row is obtained from Theorem~\ref{thm:main} and Corollaries~\ref{cor:B} and~\ref{cor:three}, the latter of which applies to $S^{n-1}$ by Remark~\ref{rem:W-sphere}; Theorem~\ref{thm:main} does not concern divergence-free vector fields.

\begin{table}[htb]
\centering
\setlength{\belowcaptionskip}{6pt}%
\caption{Numbers of generators for $S^{n-1}$, $n\geq3$.}
\label{tab:counts}
\small
\setlength{\tabcolsep}{5pt}%
\begin{tabular}{@{}l*{6}{c}@{}}

\toprule

 & \multicolumn{3}{c}{$\xalgz(S^{n-1})$ and $\operatorname{Diff}_0(S^{n-1},\omega)$} & \multicolumn{3}{c@{}}{$\xalg(S^{n-1})$ and $\operatorname{Diff}_0(S^{n-1})$}\\
\cmidrule(lr){2-4}\cmidrule(l){5-7}
 & vector fields & flows & diffeomorphisms & vector fields & flows & diffeomorphisms\\

\midrule

Lukatskii & $n$ & $3$ & $4$ & $n+1$ & $4$ & $2$\\
Lukatskii${}+{}$Kuranishi & $2$ & $2$ & $3$ & $2$ & $2$ & $3$\\
Theorem~\ref{thm:main} & --- & --- & --- & $2$ & $2$ & $3$\\
\bottomrule

\end{tabular}
\end{table}

\subsection{The conformal algebra of $S^2$ is maximal}\label{ssec:S2}

In this subsection, we prove the case $n=3$ of Theorem~\ref{thm:luk2}, i.e., that the conformal algebra $\frakc=\sX_0\oplus\{G_a\mid a\in\R^3\}$ is a maximal Lie subalgebra of $\xalg(S^2)$.
We use the notation of this section for $n=3$, in particular, the vector fields $X_f$ and $G_f$ introduced before Proposition~\ref{prop:sphere-pairs}, and further notation of~\cite{CofDS2}, which we now recall.

For $k\in\N_0$, let $H_k$ be the space of restrictions to $S^2$ of the homogeneous polynomials of degree $k$ in $x_1,x_2,x_3$ that are harmonic, i.e., annihilated by the Euclidean Laplacian $\Delta$.
Then, $\Ocal(S^2)=\bigoplus_{k\in\N_0}H_k$, $\dim H_k=2k+1$, the space $H_k$ is invariant and irreducible under $\mathrm{SO}_3(\R)$, and $\Delta_{S^2}f=-k(k+1)f$ for $f\in H_k$,
where $\Delta_{S^2}$ is the Laplace--Beltrami operator of $S^2$ (see~\cite[Section~2.1]{CofDS2} and~\cite[Theorems~1.4.5 and~1.7.2]{DaiXu}).
For $k\geq1$, let $X_{H_k}:=\{X_f\mid f\in H_k\}$, $G_{H_k}:=\{G_f\mid f\in H_k\}$, and $I_k:=X_{H_k}\oplus G_{H_k}$; the maps $f\mapsto X_f$ and $f\mapsto G_f$ are injective on $H_k$, and they commute with rotations, so each $I_k$ is invariant under $\mathrm{SO}_3(\R)$, and $[X_f,X_p]=X_{X_f(p)}$ and $[X_f,G_p]=G_{X_f(p)}$ for $f\in H_1$ and every polynomial $p$ (see~\cite[Section~2.2]{CofDS2}).
By~\cite[Proposition~1]{CofDS2}, every $V\in\xalg(S^2)$ has a unique decomposition
\begin{equation}\label{eq:S2-decomposition}
V=\sum_{k\geq1}V^{(k)}, \qquad
V^{(k)}=X_{p_k}+G_{q_k}\in I_k,
\end{equation}
with $p_k,q_k\in H_k$, and only finitely many $V^{(k)}$ are nonzero; thus, $\xalg(S^2)=\bigoplus_{k\geq1}I_k$, and $I_1=\frakc$.
Moreover, with $\operatorname{curl}V:=-\operatorname{div}(x\times V)$, the $H_k$-components of $\operatorname{curl}V$ and of $\operatorname{div}V$ are $-k(k+1)p_k$ and $-k(k+1)q_k$, respectively, since $\operatorname{div}X_f=0$, $\operatorname{div}G_f=\Delta_{S^2}f$, $\operatorname{curl}X_f=\Delta_{S^2}f$, and $\operatorname{curl}G_f=0$ for every polynomial $f$ (see~\cite[Section~2.2]{CofDS2}); here, $x\times X_f=-G_f$ and $x\times G_f=X_f$.
We write $V\mapsto V^{(k)}$ for the projection onto $I_k$ along the other summands.

\begin{lemma}\label{lem:S2-bracket}
Let $k\geq1$, $p\in H_k$, and $s_1,s_2\in\R$. Then,
$$
[G_{x_1},s_1 G_p+s_2 X_p]^{(k)}=\frac{2}{k(k+1)}\bigl(s_1 X_{p'}-s_2 G_{p'}\bigr),   \quad\mbox{where } p':=X_p(x_1)\in H_k.
$$
\end{lemma}

\begin{proof}
Let $W:=[G_{x_1},s_1 G_p+s_2 X_p]$.
Since $\operatorname{div}[U,V]=U(\operatorname{div}V)-V(\operatorname{div}U)$, $\operatorname{div}X_p=0$, $\operatorname{div}G_p=\Delta_{S^2}p=-k(k+1)p$,
and $\operatorname{div}G_{x_1}=\Delta_{S^2}x_1=-2x_1$, we have that
$$
  \operatorname{div}W=-s_1 k(k+1)\,G_{x_1}(p)+2s_1 \, G_p(x_1)+2s_2 \, X_p(x_1).
$$
On $S^2$, we have that $G_{x_1}(p)=G_p(x_1)=\partial_1 p - k x_1 p$, since $x \cdot \nabla p=kp$ and $x \cdot \nabla x_1=x_1$.
For $p\in H_k$, the polynomial $x_1 p-\frac{1}{2k+1}|x|^2\partial_1 p$ is harmonic, since $\Delta(x_1 p)=2\partial_1 p$ and $\Delta(|x|^2\partial_1 p)=(4k+2)\partial_1 p$; thus,
    on $S^2$, the function $x_1 p$ belongs to $H_{k+1}\oplus H_{k-1}$, and $\partial_1 p\in H_{k-1}$.
Hence, the $H_k$-component of $\operatorname{div}W$ is $2s_2 X_p(x_1)$, which belongs to $H_k$, since $X_p(x_1)=-X_{x_1}(p)$ and $H_k$ is invariant under rotations.
Thus, $q_k=-\frac{2s_2}{k(k+1)}X_p(x_1)$ in the decomposition~\eqref{eq:S2-decomposition} of $W$.

For the curl, note that $[G_{x_1},x\times U]=x\times[G_{x_1},U]$ for every smooth vector field $U$ on $S^2$.
To wit, for $a\in\R^3$, a direct computation with $G_a= a - \ell_a(x)x$ gives $[G_a,x\times U]-x\times[G_a,U]=(a-\ell_a(x)x)\times U + \langle a,x\times U\rangle x$, and the right-hand side vanishes on $S^2$, since $(a-\ell_a(x) x)\times U$ is normal to $S^2$ with normal component $\langle x,(a-\ell_a(x)x)\times U\rangle=-\langle a , x\times U\rangle$.
(Conceptually, the flow of $G_a$ consists of conformal maps of $S^2$, which commute with the complex structure $U\mapsto x\times U$.)
Hence, $\operatorname{curl}W=-\operatorname{div}[G_{x_1},x\times(s_1 G_p+s_2 X_p)]=-\operatorname{div}[G_{x_1},s_1 X_p-s_2 G_p]$, and the computation above, with $(s_1,s_2)$ replaced by $(-s_2,s_1)$, shows that the $H_k$-component of $\operatorname{curl}W$ is $-2s_1 X_p(x_1)$.
Thus, $p_k=\frac{2s_1}{k(k+1)}X_p(x_1)$, which completes the proof.
\end{proof}

\begin{proof}[Proof of Theorem~\ref{thm:luk2} for $n=3$]
Let $\mathfrak l$ be a Lie subalgebra of $\xalg(S^2)$ that strictly contains $\frakc$; we show that $\mathfrak l=\xalg(S^2)$ in three steps.

\xc{Step~1: $\mathfrak l=\bigoplus_{k\geq1}(\mathfrak l \cap I_k)$.}
Let $C:=\sum_{i=1}^3\operatorname{ad}_{X_{x_i}}^2$, which preserves $\mathfrak l$, since $X_{x_i} \in \frakc \subseteq \mathfrak l$.
For every polynomial $p$, the operator $\operatorname{ad}_{X_{x_i}}$ acts on $X_p$ and on $G_p$ through $p\mapsto X_{x_i}(p)$, and a direct computation gives
$$
\sum_{i=1}^3X_{x_i}^2=|x|^2\Delta-(x\cdot\nabla)^2-(x\cdot\nabla)
\quad\mbox{on }\R^3,
$$
and this operator acts on the harmonic homogeneous polynomials of degree $k$ by the scalar $-k(k+1)$.
Thus, $C$ acts on $I_k$ by the scalar $-k(k+1)$.
For $V\in\mathfrak l$, the component $V^{(k)}$ is obtained by applying to $V$ the polynomial $\prod_{j\neq k}\frac{C+j(j+1)}{j(j+1)-k(k+1)}$ in $C$, where $j$ ranges over the finitely many degrees other than $k$ for which $V^{(j)}\neq0$; hence, $V^{(k)}\in\mathfrak l$.

\xc{Step~2: if $\mathfrak l\cap I_k\neq0$ for some $k\geq1$, then $I_k\subseteq\mathfrak l$.}
Let $U:=\mathfrak l\cap I_k$, which is invariant under $\mathrm{SO}_3(\R)$: the Lie algebra $\mathfrak l$ contains $\sX_0$, so it is invariant under $\mathrm{SO}_3(\R)$ by the argument in the proof of Proposition~\ref{prop:sphere-pairs}, and so is $I_k$.
By Step~1 and $G_{x_1}\in\mathfrak l$, the map $u\mapsto[G_{x_1},u]^{(k)}$ sends $U$ into $U$.
Let $\pi_G,\pi_X\colon I_k\to H_k$ be the equivariant projections $X_p+G_q\mapsto q$ and $X_p+G_q\mapsto p$, respectively.
Let $D\colon H_k\to H_k$ be given by $Df:=X_f(x_1)=-X_{x_1}(f)$, which is, up to sign, the infinitesimal generator of the rotations about $e_1$ on $H_k$.
It is nonzero: the rotations that act trivially on $H_k$ form a normal subgroup of $\mathrm{SO}_3(\R)$; if it contained the rotations about $e_1$, it would contain all of their conjugates, i.e., all rotations, and $H_k$ would consist of constant functions, which is impossible for $k\geq1$.
Fix $p\in H_k$ with $Dp\neq0$.

If $\pi_G(U)=0$, then $U\subseteq X_{H_k}$, so $U=X_{H_k}$ by the irreducibility of $H_k$,
    and Lemma~\ref{lem:S2-bracket} with $(s_1,s_2)=(0,1)$ gives $G_{p'}\in U$ with $p'=Dp\neq0$, a contradiction.
Hence, $\pi_G(U)$ is a nonzero invariant subspace of $H_k$, so $\pi_G(U)=H_k$, and likewise $\pi_X(U)=H_k$.
If $\pi_G|_U$ has a nonzero kernel, then this kernel is a nonzero invariant subspace of $X_{H_k}$, so it is $X_{H_k}$; together with $\pi_G(U)=H_k$, this gives $U=I_k$.
Likewise, $U=I_k$ if $\pi_X|_U$ has a nonzero kernel.
Thus, we may assume that $\pi_G|_U$ and $\pi_X|_U$ are injective, hence bijective onto $H_k$.
Then, $\phi:=\pi_X\circ(\pi_G|_U)^{-1}$ is an equivariant linear automorphism of $H_k$, and
$$
U=\{G_q+X_{\phi(q)}\mid q\in H_k\}.
$$
Since $\phi$ is equivariant, it commutes with $D$.
By Lemma~\ref{lem:S2-bracket}, applied to $G_q$ and to $X_{\phi(q)}$ separately,
$$
    [G_{x_1},G_q+X_{\phi(q)}]^{(k)}=\frac{2}{k(k+1)}\bigl(X_{Dq}-G_{D\phi(q)}\bigr)\in U.
$$
Writing $X_{Dq}-G_{D\phi(q)}=G_{q'}+X_{\phi(q')}$ gives $q'=-D\phi(q)$ and $\phi(q')=Dq$, so $\phi^2Dq=-Dq$ for all $q\in H_k$.
The kernel of $\phi^2+1$ is an invariant subspace of $H_k$ that contains the nonzero image of $D$, so $\phi^2=-1$ on $H_k$ by irreducibility.
But then, $\phi$ turns $H_k$ into a complex vector space, which is impossible, since $\dim H_k=2k+1$ is odd.
This contradiction shows that $U=I_k$.

\xc{Step~3: propagation.}
Since $\mathfrak l\neq\frakc=I_1$, Steps~1 and~2 give $I_k\subseteq\mathfrak l$ for some $k\geq2$.
For $p\in H_k$, the formula for $\operatorname{div}W$ in the proof of Lemma~\ref{lem:S2-bracket}, with $(s_1,s_2)=(1,0)$, gives
$$
\operatorname{div}[G_{x_1},G_p]=\bigl(2-k(k+1)\bigr)(\partial_1p-kx_1p)
=\bigl(2-k(k+1)\bigr)\Bigl(-k\,h_{k+1}+\frac{k+1}{2k+1}\,\partial_1p\Bigr)
\quad\mbox{on } S^2,
$$
where $h_{k+1}\in H_{k+1}$ is the restriction of the harmonic polynomial $x_1p-\frac{1}{2k+1}|x|^2\partial_1p$.
The same computation appears in the proof of~\cite[Lemma~2]{Luk77}.
Now, let $p\in H_k$ be the real part of $(x_1+\mathrm{i}x_2)^k$.
Then, $p(e_1)=1$ and $\partial_1p(e_1)=k$, so that $\partial_1p\neq0$ and $h_{k+1}(e_1)=1-\frac{k}{2k+1}\neq0$.
Also, $2-k(k+1)\neq0$, since $k\geq2$.
Hence, $[G_{x_1},G_p]\in\mathfrak l$ has nonzero components in $G_{H_{k+1}}$ and in $G_{H_{k-1}}$, so $\mathfrak l\cap I_{k+1}\neq0$ and $\mathfrak l\cap I_{k-1}\neq0$ by Step~1, and Step~2 gives $I_{k+1}\subseteq\mathfrak l$ and $I_{k-1}\subseteq\mathfrak l$.
By induction, $I_j\subseteq\mathfrak l$ for all $j\geq1$, and $\mathfrak l=\xalg(S^2)$ by~\eqref{eq:S2-decomposition}.
\end{proof}

\printbibliography

\end{document}